\documentclass[11pt,a4paper,reqno]{amsart}

\makeatletter
\renewcommand{\tocsection}[3]{%
  \indentlabel{\@ifnotempty{#2}{\bfseries\ignorespaces#1 #2\quad}}\bfseries#3}
\renewcommand{\tocsubsection}[3]{%
  \indentlabel{\@ifnotempty{#2}{\ignorespaces#1 #2\quad}}#3}

\newcommand\@dotsep{4.5}
\def\@tocline#1#2#3#4#5#6#7{\relax
  \ifnum #1>\c@tocdepth 
  \else
    \par \addpenalty\@secpenalty\addvspace{#2}%
    \begingroup \hyphenpenalty\@M
    \@ifempty{#4}{%
      \@tempdima\csname r@tocindent\number#1\endcsname\relax
    }{%
      \@tempdima#4\relax
    }%
    \parindent\z@ \leftskip#3\relax \advance\leftskip\@tempdima\relax
    \rightskip\@pnumwidth plus1em \parfillskip-\@pnumwidth
    #5\leavevmode\hskip-\@tempdima{#6}\nobreak
    \leaders\hbox{$\m@th\mkern \@dotsep mu\hbox{.}\mkern \@dotsep mu$}\hfill
    \nobreak
    \hbox to\@pnumwidth{\@tocpagenum{\ifnum#1=1\bfseries\fi#7}}\par
    \nobreak
    \endgroup
  \fi}
\AtBeginDocument{%
\expandafter\renewcommand\csname r@tocindent0\endcsname{0pt}
}
\def\l@subsection{\@tocline{2}{0pt}{2.5pc}{5pc}{}}
\makeatother

\usepackage[linktocpage]{hyperref}
\hypersetup{
  colorlinks   = true, 
  urlcolor     = blue, 
  linkcolor    = Purple, 
  citecolor   = red 
}
\usepackage{changepage}
\usepackage{amsmath}
\usepackage{amsthm}
\usepackage{amsfonts}
\usepackage{amssymb}
\usepackage{array}
\usepackage[margin=2.7cm]{geometry}

\usepackage{mathtools}
\usepackage{soul}
\usepackage{enumitem}
\usepackage{hhline}
\usepackage{multirow} 

\makeatletter
\newcommand{\mylabel}[2]{#2\def\@currentlabel{#2}\label{#1}}
\makeatother

\newcommand{\C}{\mathcal{C}}
\newcommand{\A}{\mathcal{A}}

\newcommand{\mV}{\mathcal{V}}
\newcommand{\mC}{\mathcal{C}}

\newcommand{\mL}{\mathcal{L}}
\newcommand{\mH}{\mathcal{H}}
\newcommand{\mE}{\mathcal{E}}
\newcommand{\mB}{\mathcal{B}}
\newcommand{\mX}{\mathcal{X}}
\newcommand{\mD}{\mathcal{D}}

\newcommand{\F}{\mathbb{F}}

\newcommand{\E}{\mathbb{E}}
\newcommand{\K}{\mathbb K}
\newcommand{\LL}{\mathbb L}
\newcommand{\N}{\mathbb{N}}
\newcommand{\Z}{\mathbb{Z}}
\newcommand{\Fq}{\F_{q}}
\newcommand{\Fm}{\F_{q^m}}
\newcommand{\Fn}{\F_{q^n}}

\newcommand{\bs}{\boldsymbol}

\newcommand{\mat}{\K^{n \times n}}

\newcommand{\ZZ}[1]{\mathbb Z/#1\mathbb Z}

\newcommand{\dd}{\mathrm{d}}

\newcommand{\qspace}{\LL[X;\theta]/(H_\Lambda(X))}

\DeclareMathOperator{\diag}{diag}
\DeclareMathOperator{\GL}{GL}
\DeclareMathOperator{\Gal}{Gal}

\DeclareMathOperator{\rk}{rk}
\DeclareMathOperator{\srk}{srk}

\DeclareMathOperator{\Aut}{Aut}
\DeclareMathOperator{\Ann}{Ann}

\DeclareMathOperator{\End}{End}
\DeclareMathOperator{\Mat}{Mat}
\DeclareMathOperator{\wt}{wt}
\DeclareMathOperator{\Nn}{N}
\DeclareMathOperator{\gcrd}{gcrd}
\DeclareMathOperator{\lclm}{lclm}
\DeclareMathOperator{\ch}{char}
\DeclareMathOperator{\ev}{ev}
\DeclareMathOperator{\Ext}{Ext}
\DeclareMathOperator{\Tr}{Tr}
\DeclareMathOperator{\vv}{v}
\DeclareMathOperator{\MM}{m}
\DeclareMathOperator{\HH}{H}
\DeclareMathOperator{\rr}{r}

\newcommand{\evB}{\ev_{\bs\mB}}

\newcommand{\Trace}{\Tr_{\LL/\K}}

\newcommand{\NN}[2]{\Nn_{#1}^{#2}}
\newcommand{\Pro}[2]{\mathcal P_{#1}[#2]}
\newcommand{\Norm}{\Nn_{\LL/\K}}
\newcommand{\xt}{\xi_\theta}

\theoremstyle{definition}
\newtheorem{theorem}{Theorem}[section]

\newtheorem{proposition}[theorem]{Proposition}
\newtheorem{corollary}[theorem]{Corollary}
\newtheorem{lemma}[theorem]{Lemma}

\newtheorem{definition}[theorem]{Definition}
\newtheorem{example}[theorem]{Example}
\newtheorem{remark}[theorem]{Remark}

 \usepackage[dvipsnames,table]{xcolor}
 \usepackage{pagecolor}
 \definecolor{light-gray}{gray}{0.90}

\usepackage{tikz}
\usetikzlibrary{arrows.meta}
\usepackage{lscape}
\usepackage{arydshln}

\newcommand\red[1]{{\color{red} #1 }}

\title[Twisted Linearized Reed-Solomon Codes]{Twisted Linearized Reed-Solomon Codes: \\A Skew Polynomial Framework}

\author[A. Neri]{Alessandro Neri}
\address{Max-Planck-Institute for Mathematics in the Sciences, \textnormal{Inselstraße 22, 04103 Leipzig, Germany}}
\email{alessandro.neri@mis.mpg.de}



\begin{document}

\maketitle
\thispagestyle{empty}
\begin{abstract}
    We provide an algebraic description for sum-rank metric codes, as quotient space of a skew polynomial ring. This approach generalizes at the same time the skew group algebra setting for rank-metric codes and the polynomial setting for codes in the Hamming metric. This allows  to construct twisted linearized Reed-Solomon codes, a new family of maximum sum-rank distance codes extending at the same time Sheekey's twisted Gabidulin codes in the rank metric and twisted Reed-Solomon codes  in the Hamming metric.  Furthermore, we provide an analogue in the sum-rank metric of Trombetti-Zhou construction, which also provides a family of maximum sum-rank distance codes. As a byproduct, in the extremal case of the Hamming metric, we obtain a new family of additive MDS codes over quadratic fields.
\end{abstract}

\tableofcontents

 \section{Introduction}
 The sum-rank metric has emerged  recently as a solution to improve the performance of  multishot network coding based on rank-metric codes \cite{nobrega2010multishot}. However, traces of this metric can be already found in connection with space-time coding \cite{el2003design,lu2005unified}. In addition, sum-rank metric codes have been shown to have applications also in the context of distributed storage, for constructing partial MDS codes \cite{martinez2019universal,martinez2020general}. Furthermore, the sum-rank distance has also been used as a metric on convolutional codes for network streaming \cite{mahmood2016convolutional}. 
 
In the last few years, the theory of codes in the sum-rank metric has been developing thanks to the works of Mart\'inez-Pe\~nas \cite{martinez2018skew,martinez2020hamming,martinez2020general}.  These codes can be seen as a generalization of Hamming-metric codes and rank-metric codes. The sum-rank metric is defined on vectors over a (non necessarily finite) field, by partitioning the ambient vector space in several blocks, considering in each of them the rank metric over a subfield and then  summing up the distances of each block. In this setting, in the case of a single block one obtains the rank metric, while partitioning the ambient space in blocks of length one, one gets the Hamming metric. In a slightly different  setting, the sum-rank metric can be considered on tuples of matrices, as it is done e.g. in \cite{byrne2020fundamental}. This setting is slightly more general, although it is equivalent when working over finite fields; see e.g. Figure \ref{fig:1}. 

Codes with optimal parameters for the sum-rank metric are called maximum sum-rank distance (MSRD for short) codes and they are essentially the analogue of MDS codes for the Hamming metric and of MRD codes for the rank metric. Only a few constructions of MSRD codes are known so far. The most prominent one is undoubtedly the family of \emph{linearized Reed-Solomon codes}, introduced in \cite{martinez2018skew}. These codes are a sort of ``tensorization'' between Reed-Solomon codes in the Hamming metric and Gabidulin codes in the rank metric. Recently, a new construction has been proposed in \cite{martinez2020general}, which generalizes linearized Reed-Solomon codes. 

\medskip

In this paper we introduce the family of \emph{twisted linearized Reed-Solomon codes}, and show that these codes are MSRD codes; see Theorem \ref{thm:twisted_are_MSR}. This family represents the natural counterpart in the sum-rank metric of \emph{twisted Gabidulin codes}  introduced by Sheekey in \cite{sheekey2016new} and of \emph{twisted Reed-Solomon codes}, defined later by Beelen, Puchinger and Rosenkilde in \cite{beelen2017twisted}. Indeed, when we specialize our construction to only one block, we obtain exactly twisted Gabidulin codes in the rank metric, while if we choose all the blocks to have length one, our codes turn out to be twisted Reed-Solomon codes in the Hamming-metric. 

We then  generalize  Trombetti and Zhou construction of maximum rank distance codes to a new family of sum-rank metric codes, called \emph{twisted linearized Reed-Solomon codes of TZ-type}. Also in this case we prove in Theorem \ref{thm:TZ_are_MSR} that the resulting codes are MSRD codes. In the special case of block length equal to one, the construction can be adapted to obtain a new family of additive MDS codes in the Hamming metric. This result is certainly of independent interest; see Theorem \ref{thm:TZ_MDS}. These codes are defined over a quadratic finite field of order $q^2$, but they are only linear over the subfield $\Fq$. Their length can be extended up to $\frac{q^2-1}{2}$.

Our constructions strongly rely on a new framework that we develop for studying codes in the sum-rank metric. We indeed propose to study codes in the sum-rank metric as subsets of a particular quotient algebra of skew polynomials. In Theorem \ref{thm:isomorphism_skew_sumrank} we show that this setting is isometric to the natural frameworks of sum-rank metric codes. Our setting unifies the linearized polynomial framework used to construct rank-metric codes over finite fields (see e.g. \cite{sheekey2016new}) -- and of the more general representation as skew polynomials or as elements of a skew algebra  used for rank-metric codes over arbitrary fields (see \cite{augot2013rank,augot2020rank}) -- with the polynomial framework of linear codes in the Hamming metric.  We remark that a similar point of view was used by Sheekey in \cite{sheekey2020new}, in which he used this approach to construct new optimal codes in the rank metric over field extensions. 

The key idea on which our framework is based, is due to recent results characterizing the roots of linearized and projective polynomials appeared in \cite{caruso2019residues,csajbok2019characterization,mcguire2019characterization}. We generalize their results over arbitrary fields, making more explicit  the connection between the sum-rank metric and the eigenspaces of a special linear map; see Theorem \ref{thm:bound_degree_improved}. Thanks to these results, we show that our framework is  isometric to the classical ones, where sum-rank metric codes are represented as sets of vectors, or as sets of tuples of matrices. We also explicitly point out how to go from one setting to another one and viceversa. Furthermore, we show how the natural nondegenerate bilinear map on the skew polynomial setting simplifies the study of duality of codes in the sum-rank metric. 

\medskip

The paper is structured as follows. Section \ref{sec:preliminaries} contains the basic notions and the mathematical background needed for the purpose of the paper. We first recall the notion of skew algebras obtained from a Galois extension $\LL/\K$ and relate them with the ring of $\K$-linear endomorphisms of $\LL$ and with the ring of the skew polynomials. In these frameworks, we give the basic definitions of rank and sum-rank metric. In Section \ref{sec:projective_polynomials} we study projective $\theta$-polynomials over arbitrary fields. We show their main properties and how one can get from them  information about the kernel of some related endomorphisms. Section \ref{sec:skewpoly_framework} introduces the new skew polynomial framework for sum-rank metric codes. This is obtained as a quotient of a skew polynomial ring, and it is shown to be isometric to the space of $\ell$-uples of square matrices. Furthermore, we define sum-rank metric codes  in this setting, where linearized Reed-Solomon codes can be seen as the natural analogue of Reed-Solomon and Gabidulin codes, generated by all the skew polynomials of bounded degree. Here we also  provide the study of equivalence of sum-rank metric codes  and we investigate the notions of adjoint of a skew polynomial and its associated Dickson matrix. In Section \ref{sec:duality} we deepen the notion of duality of sum-rank metric codes. Our definition is shown to be equivalent to the natural dualities defined on the space of vectors and on the space of $\ell$-uples of matrices. Furthermore, the skew polynomial framework simplifies drastically the study of the duals of linearized Reed-Solomon codes. We then introduce the family of {twisted linearized Reed-Solomon codes} in Section \ref{sec:twisted_LRS}. We show that they constitute a family of MSRD codes and then derive their dual and adjoint codes. In Section \ref{sec:TZ} we investigate the family twisted linearized Reed-Solomon codes of TZ-type. After showing that they are MSRD, we focus on the special case of blocks of length one, and the new obtained MDS codes in the Hamming metric.  Finally, we draw our conclusions and list some open problems in Section \ref{sec:conclusions}.

\section{Preliminaries}\label{sec:preliminaries}

\subsection{Skew Group Algebras}\label{sec:skew_group_algebra}

Let $\LL$ be a field, let  $G$ be a finite group and let $\phi:G\rightarrow \Aut(\LL)$ be a group homomorphism. We denote by $(\LL^{\phi}[G],+,\circ_\phi)$ the skew group algebra
$$\LL^{\phi}[G]\coloneqq\bigg\{\sum_{g\in G}a_gg \mid a_g \in \LL \bigg\},$$
with the usual addition $+$ that works componentwise, i.e.
$$\sum_{g \in G}a_gg + \sum_{g\in G}b_gg=\sum_{g \in G}(a_g+b_g)g,$$
and where the multiplication $\circ_{\phi}$ is defined as
$$a_g g \circ_{\phi} b_hh\coloneqq a_g(\phi(g)(b_h))(gh),$$ and then extended by associativity and distributivity.
Observe that if we take $\phi$ to be the trivial homomorphism mapping each element of $G$ to the identity automorphism of $\LL$, we obtain the classical group algebra $\LL G$. 

We now fix the following framework. Let $\LL/\K$ be  a Galois extension whose Galois group is $G\coloneqq \Gal(\LL/\K)$. Since $G\leq \Aut(\LL)$, we consider the natural inclusion map $ G\hookrightarrow \Aut(\LL)$ and we denote it by $\iota$. Clearly, $\iota$ is a group homomorphism and one can clearly consider the skew group algebra $\LL[G]\coloneqq \LL^\iota[G]$. 

The skew group algebra $\LL[G]$ has been recently studied in connection with rank-metric codes \cite{augot2020rank,elmaazouz2021}. Recall that the theory of rank-metric codes deals with the study of subspaces of the matrix space $\K^{n\times n}$ equipped with the rank metric, where the distance between two matrices $A,B \in \K^{n\times n}$ is  $\rk(A-B)$.
The reasons why $\LL[G]$  is fundamental in the study of these codes are the following.
First, it is isomorphic to the $n\times n$ matrix algebra over $\K$. More precisely, for any element $a=\sum_g a_gg \in \LL[G]$, one can associate a map 
$$\begin{array}{rcl}
\psi_a:\LL &\rightarrow &\LL \\
   \beta & \longmapsto & \sum_g a_gg(\beta).
\end{array}$$

\begin{theorem}[\textnormal{\cite[Theorem 1.3]{chase1969galois}}]\label{thm:Galois_endomorphism}
 The map $:a\longmapsto \psi_a$ is a $\K$-algebra isomorphism between $\LL[G]$ and $\End_\K(\LL)$. 
\end{theorem}

With a slight abuse of notation, from now on we will write $a$ also to indicate the map $\psi_a$, and so we will refer to $\rk(a)$ and $\ker(a)$ to indicate the rank over $\K$ of $\psi_a$ and its kernel, respectively.

The second reason why this skew group algebra approach is useful for studying rank-metric codes is described in the following. In the theory of codes equipped with the rank metric, the classical framework is given by choosing as fields $\K=\Fq$ and $\LL=\Fm$. With these assumptions, $G=\Gal(\Fm/\Fq)$ is always cyclic and when this happens,   the skew group algebra $\LL[G]=\LL[\theta]$, where $\theta$ is any generator of $\Gal(\LL/\K)$. Hence, the elements of $\LL[\theta]$ are represented as polynomials in $\theta$ with coefficients in $\LL$ and are also called \textbf{$\theta$-polynomials} (or \textbf{$\theta$-linearized polynomials}). This naturally induces a notion of $\theta$-degree of a nonzero $\theta$-polynomial $f\coloneqq \sum_i f_i \theta^i$, as $\deg_\theta(f)\coloneqq \max\{i \mid f_i \neq 0\}$. This notion allows to establish a fundamental result that gives an upper bound on the nullity of a $\theta$-polynomial.

\begin{theorem}[\textnormal{\cite{ gow2009galoisext,augot2013rank}}]\label{thm:fundamental_theorem_theta_poly}
 Let $\LL/\K$ be a cyclic Galois extension whose Galois group is generated by $\theta$, and let $f\in \LL[\theta]$ be a nonzero $\theta$-polynomial. Then
 $$\dim_\K(\ker(f))\leq \deg_\theta(f).$$
\end{theorem}

Theorem \ref{thm:fundamental_theorem_theta_poly} is the analogue of the well-known result on polynomials over a field that states that the number of roots of a nonzero polynomial is at most its degree. 
A generalization of this result to  polynomials in $m$ variables has been provided in \cite{augot2020rank}, where instead of having a cyclic extension, one considers an Abelian extension whose Galois group is isomorphic to $\ZZ{n_1}\times \ldots \times \ZZ{n_m}$. However, we will later improve this result in Theorem \ref{thm:bound_degree_improved}, give a more refined bound.

Lastly, representing the space $\K^{n \times n}$ via the skew group algebra $\LL[\theta]$ allows to give a natural notion of $\LL$-linearity to rank-metric codes. In the case of finite fields, it was shown that rank-metric codes with optimal parameters (called \textbf{maximum rank distance codes}) are dense among $\LL$-linear codes \cite{neri2018genericity}. 

We conclude this section by recalling the notion of {adjoint} of a $\theta$-polynomial. For a given $f =f_0\mathrm{id}+\ldots+f_{n-1}\theta^{n-1}\in \LL[\theta]$, the \textbf{adjoint} of $f$ is the $\theta$-polynomial
\begin{equation}\label{eq:adjoint} f^\top\coloneqq \sum_{i=0}^{n-1}{\theta^{n-i}}(f_i) \theta^{n-i}\in \LL[\theta].\end{equation}
Note that the term ``adjoint'' comes  from the fact that it represents the adjoint operator on $\LL$ with respect to the trace bilinear form, that is,
 $$ \Trace(f(\alpha)\beta)=\Trace(\alpha f^\top(\beta)),\qquad \mbox{ for every } \alpha, \beta \in \LL.$$
In the identity above, the map $\Trace$ is the \textbf{trace map} with respect to the extension $\LL/\K$, defined as
$$\Trace(\alpha)\coloneqq \sum_{\sigma \in \Gal(\LL/\K)} \sigma(\alpha)=\sum_{i=0}^{n-1}{\theta}^i(\alpha),$$  

\medskip

\subsection{Norm of a Galois Extension}

In this subsection we study the norm associated to a field extension $\LL/\K$. It can be definied over any arbitrary extension, but for our purposes we only consider when $\LL/\K$ being Galois. In this case, the norm of an element is given by the product of all its conjugates under the Galois group action, and it has well-known and interesting properties that we summarize here.

\begin{definition}
  Let $\LL/\K$ be a cyclic Galois extension of 
  field and let $\theta$ be a generator of $\Gal(\LL/\K)$. For $\alpha \in \LL$, the \textbf{norm} of $\alpha$ with respect to the extension $\LL/\K$ is defined by
  $$\Norm(\alpha)  \coloneqq  \prod_{\sigma \in \Gal(\LL/\K)} \sigma(\alpha)=\prod_{i=0}^{n-1}{\theta}^i(\alpha),$$
\end{definition}
We will refer to the function 
$$\begin{array}{rcl} \Norm:\LL & \longrightarrow & \K
\end{array}$$
as the \textbf{norm map} of $\LL/\K$. Moreover, for any nonnegative integer $i$, define the \textbf{$i$-th truncated norm} as the  map 
$$ \begin{array}{rcl} \NN{\theta}{i}:\LL & \longrightarrow & \LL \\
\beta & \longmapsto & \prod\limits_{j=0}^{i-1} \theta^j(\beta).\end{array}$$

Now, by a well-known result of Hilbert, we introduce the last map, which connects the norm and the $i$-th truncated norm. For any $\theta$ generator of $\Gal(\LL/\K)$, we define 

$$\begin{array}{rcl}
\xt :\LL^* & \longrightarrow &\LL^* \\
 \alpha & \longmapsto & \frac{\theta(\alpha)}{\alpha}.
\end{array}$$

The following results are well known and their proofs can can be found in any Algebra textbook.

\begin{lemma}\label{lem:normproperties}
   The norm map of $\LL/\K$ satisfies the following properties:
   \begin{enumerate}
   \item $\Norm(\alpha) \in \K$ for all
     $\alpha \in \LL$.
 \item $\Norm(\alpha)=0$ if and only if $\alpha=0$.
   \item $\Norm$ restricted to $\LL^*$ is a surjective group homomorphism from $\LL^*$ to $\K^*$.
 \item $\xt$ is a group homomorphism from $\LL^*$ to
   itself.
 \item For every generator $\theta$ of $\Gal(\LL/\K)$, $\xt(\alpha)=1$ if and only
   if $\alpha \in \K^*$.
 \item (Hilbert's Theorem 90) For every generator $\theta$ of $\Gal(\LL/\K)$, $\ker (\Norm)=\mathrm{Im}(\xt)$.
 \item For every generator $\theta$ of $\Gal(\LL/\K)$ and for every $\K$-subspace $V$ of $\LL$, 
 $$ \xt(V\setminus\{0\})\cong (V\setminus\{0\})/\K^*.$$
 \end{enumerate}
\end{lemma}

\begin{lemma}\label{lem:normpreimage}
 Let $\alpha \in \LL^*$ and $\theta$ be a generator of $\Gal(\LL/\K)$. Then 
\begin{enumerate}
\item If $\alpha \notin \ker(\Norm)$, then $\xt^{-1}(\{\alpha\})=\emptyset$.
\item Let $\alpha \in \ker(\Norm)$. If $x_1, x_2 \in \xt^{-1}(\{\alpha\})$, then $\frac{x_1}{x_2} \in \K^*$, or equivalently, there exists an $x \in \Fm$ such that
$$\xt^{-1}(\{\alpha\})=\left\{ \lambda x\mid \lambda \in \K^* \right\}.$$

Moreover such an $x$ is of the form
$ x= \chi_\alpha(\gamma^{-1}),$
where 
$$\chi_\alpha=\sum_{i=0}^{n-1}\NN{\theta}{i}(\alpha)\theta^{i},$$
and $\gamma \in \LL^*$ is such that $\chi_\alpha(\gamma)\neq 0$.
\end{enumerate}
\end{lemma}

 Note that in the statement of Lemma \ref{lem:normpreimage}(2)  there always exists an element $\gamma
\in \LL^*$ such that $\chi_\alpha(\gamma)\neq 0$. This is due to Theorem \ref{thm:Galois_endomorphism}, since $\chi_\alpha \in \LL[\theta]\setminus\{0\}$.

\medskip

\subsection{Skew Polynomials}

In this section we introduce the rings of skew polynomials and recall some of their basic properties. These polynomials naturally arises studying the skew group algebra $\LL[\theta]$, which can be obtained as a quotient of this ring. There is plenty of literature on skew polynomial rings, starting from the pioneering work of Ore \cite{ore1933theory}. 

Let $\LL$ be a field and let $\theta \in \Aut(\LL)$. Define the skew polynomial ring $(\LL[X;\theta],+,\cdot)$ as the set of polynomials in the indeterminate $X$ with coefficients in $\LL$, where ``$+$'' is the usual polynomial addition, while the multiplication  follows the rule
$$ X a=\theta(a) X, \mbox{ for any } a \in \LL,$$
and it is then extended by associativity and distributivity.

It is well-known that the ring $\LL[X;\theta]$ is commutative if and only if $\theta=\mathrm{id}$, and in general its center is $\K[X^m]$, where $\K\coloneqq \LL^{\theta}$ is the subfield of $\LL$ fixed by $\theta$ and $m$ is the order of $\theta$.  
Moreover, 
 the ring $\LL[X;\theta]$ is a right-Euclidean domain, and hence there is a well-defined notion of greatest common right divisor and least common left multiple between two skew polynomials $F_1$ and $F_2$, which are denoted, respectively, by $\gcrd(F_1, F_2)$ and $\lclm(F_1,F_2)$. Throughout the paper, the right-divisibility of polynomials will be written as $F(X)\mid_{\rr} G(X)$.

From now on we will focus on a field $\LL$ with an automorphism $\theta$ of order $n$, and call $\K\coloneqq \LL^\theta$. Thus, we are in the same setting as the one introduced in Section \ref{sec:skew_group_algebra}, where $\LL/\K$ is a cyclic Galois extension, and $\theta$ is a generator of $\Gal(\LL/\K)$. In this setting, define
  the map
  \begin{equation}\label{eq:psi_skewpoly-skewalgebra}
      \begin{array}{rccc}
       \Phi : &\LL[X;\theta] & \longrightarrow & \LL[\theta] \\
       &f_0+f_1X+\ldots +f_dX^d & \longmapsto & f_0\mathrm{id}+f_1\theta+\ldots+f_d\theta^d.
      \end{array}
  \end{equation}

The importance of the map $\Phi$ is given by the following well-known result.

\begin{theorem}\label{thm:ismorphism_skew_theta}
 The map $\Phi$ defined in \eqref{eq:psi_skewpoly-skewalgebra}
 is a $\K$-algebra surjective homomorphism (and an $\LL$-vector space homomorphism) whose kernel is the twosided ideal 
  $(X^n-1)$. Consequently, $$ \LL[X;\theta]/(X^n-1) \cong \LL[\theta].$$ 
\end{theorem}

\noindent \textbf{Notation.} As a consequence of Theorem \ref{thm:ismorphism_skew_theta}, from now on we will abbreviate the notation, writing $\ker(F)\coloneqq \ker(\Phi(F))$, for any skew polynomial $F\in\LL[X;\theta]$. Moreover, we can define the \textbf{evaluation} of a skew polynomial $F$ in an element $\beta\in \LL$ as the evaluation of $\Phi(F)$ in $\beta$, and write $F(\beta)$.

\begin{corollary}\label{cor:kerf_skew}
 Let $F\in \LL[X;\theta]$. Then
 $$ \dim_\K(\ker (F))=\deg(\gcrd(F(X),X^n-1)).$$
\end{corollary}

By Corollary \ref{cor:kerf_skew} we can see that the degree of the greatest common right divisor of a skew polynomial $F(X)$ with $X^n-1$ reveals important information about the $\K$-linear  map induced by $F(X)$. We can generalize this approach, and substitute $X^n-1$ with any other central skew polynomial $X^n- \lambda$, for $\lambda \in\K$.

\begin{definition}
 Let $F\in \LL[X;\theta]$, and let $\lambda \in \K^*$. We define the \textbf{$\lambda$-value for $F$} to be the integer 
 $$ d_\lambda(F)\coloneqq \deg(\gcrd(F(X), X^n-\lambda)).$$
\end{definition}

One can immediately see that by Corollary \ref{cor:kerf_skew} we have $d_1(F)= \dim_{\K}(\ker(F))$. However, all the $\lambda$-values represent important quantities for the study of the map given by $F$, and especially for the corresponding $\theta$-projective polynomial that we introduce in the next section.

\medskip

\subsection{Rank and Sum-Rank Metric}

In the last decade rank-metric codes have  gained attention due to their numerous applications, such as in random network coding \cite{silva2008rank}, cryptography \cite{gabidulin1991ideals} and distributed storage \cite{rawat2013optimal}. However, codes in the rank metric have been already introduced  by Delsarte \cite{delsarte1978bilinear}, Gabidulin \cite{gabidulin1985theory} and Roth \cite{roth1991maximum} independently. Originally introduced as sets of $n\times m$ matrices over a finite field $\Fq$, they can be equivalently represented as sets of vectors of length $n$ over a field extension $\Fm$. However, when we restrict to $n=m$, another equivalent representation is via the ring of $q$-linearized polynomials with coefficients in $\Fn$ modulo the ideal generated by $x^{q^n}-x$; see \cite{sheekeysurvey} for a survey on rank-metric codes from this viewpoint. The representation of rank-metric codes as $q$-linearized polynomials is exactly the same as the one via the skew algebra $\Fn[\theta]$, as explained in Section \ref{sec:skew_group_algebra}. We briefly describe rank-metric codes and consequently sum-rank metric codes in the setting of skew algebras.

Let $\theta$ be a generator of $\Gal(\LL/\K)$. The \textbf{rank metric} on $\LL[\theta]$ is the distance defined by the rank, that is
$$ \dd_{\rk}(f,g):=\rk(f-g), \qquad \mbox{ for ever } f,g, \in \LL[\theta].$$
A \textbf{rank-metric code} $\C$ is a subset of $\LL[\theta]$ endowed with the rank metric. The \textbf{minimum rank distance} of $\mC$ is the integer
$$ \dd_{\rk}(\C):=\min\{\dd_{\rk}(f,g) : f,g \in \C, f\neq g \}.$$
This space is often considered endowed with a nondegenerate bilinear form. We denote by $\langle \cdot , \cdot \rangle_{\rk}$ the nondegenerate $\K$-bilinear form on $\LL[\theta]$, given by
$$ \langle a, b\rangle_{\rk}\coloneqq\Trace\bigg(\sum_{i=0}^{ n-1} a_ib_i\bigg)=\Trace\left(\langle a, b\rangle_{\rk,\LL}\right), $$
for any $a=a_0\mathrm{id}+\ldots+a_{n-1}\theta^{n-1}$, $b=b_0\mathrm{id}+\ldots+b_{n-1}\theta^{n-1} \in \LL[\theta]$. This bilinear form is often used to define a notion of \emph{dual codes}. The interested reader is referred to \cite{augot2020rank} for a comprehensive study of the skew algebra setting for rank-metric codes.
\medskip

This setting can be naturally generalized  in order to define the \textbf{sum-rank metric}. We first introduce the classical frameworks in which sum-rank metric codes are usually considered.
The first case of codes interpreted as vectors was deeply analyzed by Mart{\'\i}nez-Pe{\~n}as \textit{et al.} in several papers; see e.g. \cite{martinez2018skew, martinez2019universal}.  The sum-rank weight on $\LL^{\ell n}$ with respect to the partition $\ell n=n+n+\ldots+n$ is defined as
$$\wt_{\vv}((v^{(1)} \mid \cdots \mid v^{(\ell)}))=\sum_{i=1}^\ell \rk_{\K}(v^{(i)}), \qquad \mbox{ for } v^{(1)},\ldots, v^{(\ell)}\in\LL^n,$$
where $\rk_{\K}(u)=\dim_\K\langle u_1,\ldots,u_n\rangle_{\K}$, for any $u=(u_1,\ldots,u_n)\in \LL^n$.
A \textbf{(vector) sum-rank code} is a subset of $\LL^{\ell n}$ endowed with the sum-rank distance $\dd_{\vv}$ induced by $\wt_{\vv}$.

As a natural generalization, one can also think as sum-rank metric codes as sets made of $\ell$-uples of matrices. More specifically, the sum-rank metric on $(\mat)^\ell$ is defined as
$$ \dd_{\MM}((A_1,\ldots,A_\ell),(B_1,\ldots, B_\ell))\coloneqq \sum_{i=1}^\ell\rk(A_i-B_i), \quad \mbox{ for  } A_1,\ldots,A_\ell,B_1,\ldots, B_\ell \in \mat.$$
Here, a \textbf{(matrix) sum-rank metric code} is a subset of $(\mat)^\ell$ endowed with the sum-rank distance $\dd_{\MM}$. A deep study of the  properties of codes in this framework has been carried out in \cite{byrne2020fundamental,gluuesing2020anticodes}.

Now, we introduce the skew algebra point of view. Consider $\ell$ copies of the space $\LL[\theta]$ and define the sum-rank distance on $(\LL[\theta])^\ell$ as
$$ \dd_{\srk}((f_1,\ldots,f_\ell),(g_1,\ldots,g_\ell))=\sum_{i=1}^\ell\dd_{\rk}(f_i,g_i).$$
Here, the duality is given by the bilinear form 
$$ \langle \bs f, \bs g\rangle_{\srk}\coloneqq\Trace\bigg(\sum_{i=1}^{\ell}\langle f_i,g_i\rangle_{\rk,\LL}\bigg)=\Trace\left(\langle \bs f, \bs g\rangle_{\srk,\LL}\right), $$
for any $\bs f=(f_1,\ldots,f_\ell)$, $\bs g=(g_1,\ldots,g_\ell) \in (\LL[\theta])^\ell$. Due to the isometry between $(\LL[\theta],\dd_{\rk})$ and $(\K^{n \times n},\rk)$ discussed in Section \ref{sec:skew_group_algebra}, we obtain that $((\LL[\theta])^{\ell},\dd_{\srk})$ is isometric to the space $(\K^{n \times n})^\ell$ endowed with the sum-rank metric defined there (see also Section \ref{sec:comparison_frameworks}). 

However, our aim is to study another framework for codes in the sum-rank metric, which will be fully described in Section \ref{sec:skewpoly_framework}.

\section{$\theta$-Projective Polynomials}\label{sec:projective_polynomials}

Analogously to the case of finite field, we consider the \emph{projective $\theta$-polynomials}, extending the definition to any cyclic Galois extension $\LL/\K$. Projective polynomials were introduced and studeid by Abhyankar in \cite{abhyankar1997projective} as polynomials over function fields in positive characteristic. Over finite fields, the number of zeros of special projective polynomials  was  investigated by many authors;\cite{bluher2004xq+,kim2021solving}. More recently, McGuire and Sheekey in \cite{mcguire2019characterization} and Csajb{\'o}k, Marino,  Polverino  and Zullo in \cite{csajbok2019characterization}, gave a constructive criterion to fully determine the number of roots of any projective polynomial over a finite field. 

\medskip

\subsection{Roots of $\theta$-Projective Polynomials}
We give an alternative point of view here which allows to generalize their notion to fields of any characteristic, and preserves many of the properties of the original projective polynomials. However, the objects we define will no longer be proper polynomials.

\begin{definition}
 A \textbf{$\theta$-projective polynomial} over $\LL$  is a map $P:\LL\rightarrow \LL$ of the form
 $$ P\coloneqq \sum_{i=0}^d a_i\NN{\theta}{i}, \quad a_i \in \LL,$$
 acting as $\beta \longmapsto \sum_{i=0}^d a_i\NN{\theta}{i}(\beta)$.  The space of $\theta$-projective polynomials over $\LL$ is denoted by $\Pro{\LL}{\theta}$. 
\end{definition}

\begin{definition}
 Let $P\in \mathcal P_\LL[\theta]$ be a $\theta$-projective polynomial. An element $\beta \in \overline{\LL}$ is said to be a \textbf{root} of $P$ if $P(\beta)=0$. Moreover, we denote the set of roots of $P$ lying in $\LL$ by $\mV_{\LL}(P)$
\end{definition}

It is natural to define the following map, which relates the skew polynomial ring $\LL[X;\theta]$ and the space of $\theta$-projective polynomials $\Pro{\LL}{\theta}$:
$$\begin{array}{rccc}
\Psi : &\LL[X;\theta] & \longrightarrow & \Pro{\LL}{\theta} \\
&\sum\limits_{i=0}^d f_i  X^i & \longmapsto & \sum\limits_{i=0}^d f_i  \NN{\theta}{i}.
\end{array}$$
In order to lighten the notation, from now on we will always write $P_F$ to denote the projective polynomial $\Psi(F)$. The connection between a skew polynomial $F$ and the corresponding $\theta$-projective polynomial $P_F$ has been already investigated in \cite{chaussade2009skew,csajbok2019characterization,mcguire2019characterization} for finite fields. 

In the following we go through many of the properties of $\theta$-projective polynomials over $\LL$, linking them to the corresponding $\theta$-polynomial, similarly as done in \cite{mcguire2019characterization} for finite fields. In particular, one can easily prove the following result relating the evaluation of a skew polynomial $F(X)$ with the evaluation of  $P_F$.

\begin{proposition}\label{prop:theta-projective_poly}
Let $F(X) \in \LL[X;\theta]$ be a nonzero skew polynomial. Then, for every $\beta \in \LL$ we have \begin{equation}\label{eq:theta-projective_roots}
F(\beta)=\beta\cdot(P_F\circ \xt (\beta) ).
\end{equation}
As a consequence, we have
$$ \ker(F)=\xt^{-1}(\mV_{\LL}(P_F))\cup\{0\} $$
\end{proposition}

\begin{proof}
 We can verify \eqref{eq:theta-projective_roots} simply by evaluating the right-hand side, since $\NN{\theta}{i}(\xt(\beta))=\theta^{i}(\beta)\beta^{-1}$.  The second equality directly follows.
\end{proof}

Now, given a skew-polynomial $F(X) \in \LL[X;\theta]$ and an element $\alpha \in \LL^*$, denote by $F_\alpha(X)$ the skew polynomial $F(\alpha X)$. One can directly verify that, if $F(x)=f_0+f_1X+\ldots+f_dX^d$, then  
\begin{equation}\label{eq:skew_Falpha} F_\alpha(X)\coloneqq  \sum_{i=0}^{d} f_i \NN{\theta}{i}(\alpha)X^i. \end{equation}

\begin{remark}
 It is easy to see that for every skew polynomial $F(X)\in \LL[X;\theta]$ and every $\alpha \in \LL$,  the map $F(X) \longmapsto F(\alpha X)$  is a ring homomorphism. Therefore, we have 
 \begin{equation}\label{eq:skew_homomoprhims_alpha}
     F_\alpha G_\alpha=(FG)_\alpha.
 \end{equation}
\end{remark}

 \begin{proposition}\label{prop:projective_shift}
 Let $F \in \LL[X;\theta]$, and $\alpha \in \LL^*$. Then, it holds that $
 P_{F_\alpha}=P_F\circ \MM_\alpha,$
where $\MM_\alpha$ denotes the map given by the multiplication by $\alpha$.
In particular, 
 $$\mV_{\LL} (P_{F_\alpha})=\alpha^{-1}\cdot \mV_{\LL}(P_F).$$
 and
 $$\ker(F_\alpha)
 =\xt^{-1}(\mV_{\LL}(P_{F_\alpha}))= \xt^{-1}(\alpha^{-1}\cdot \mV_{\LL}(P_F))\cup\{0\}.$$
\end{proposition}

\begin{proof}
 This is an easy calculation. By \eqref{eq:skew_Falpha} one can immediately see that for every $\beta\in\LL$ we have
 $$ P_{F_\alpha}(\beta)=\sum_{i=0}^{d} f_i \NN{\theta}{i}(\alpha)\NN{\theta}{i}(\beta)=\sum_{i=0}^{d} f_i \NN{\theta}{i}(\alpha\beta)=P_F(\alpha\beta).$$
 The other two identities follow immediately.
\end{proof}

{The following result clarify why the term \emph{$\theta$-projective polynomial} is indeed appropriate.}

\begin{theorem}
Let $F=f_0+\ldots+f_d X^{d}$ be a nonzero $\theta$-polynomial with $f_0\neq 0$. Then,
 $$\mV_{\LL}(P_F)=\bigsqcup_{\lambda \in \K^*} \mX_\lambda,$$
 where 
 $\mX_\lambda=\alpha_\lambda \cdot \xt(\ker(F_{\alpha_{\lambda}})\setminus\{0\})$ is a copy of $\mathbb P^{d_{\lambda}(F)-1}(\K)$, and  $\alpha_\lambda\in \LL^*$ is any element such that $\Norm(\alpha_\lambda)=\lambda$.
\end{theorem}

\begin{proof}
 For each $\lambda \in\K^*$, let us fix a representative of the norm function in $\LL^*$, i.e. an element $\alpha_{\lambda}\in \LL^*$ such that $\Norm(\alpha_\lambda)=\lambda$. Let $\beta\in \LL^*$ such that $\Norm(\beta)=\lambda$. This means that there exists $\mu\in \LL^*$ such that $\beta=\alpha_\lambda \xt(\mu)$, and such an element $\mu$ is unique modulo scalar multiplication in $\K^*$. Now, observe that
 $$P_F(\beta)=P_F(\alpha_\lambda \xt(\mu))\stackrel{(*)}{=}P_{F_{\alpha_\lambda}}(\xt(\mu))\stackrel{(**)}{=}\mu^{-1}F_{\alpha_{\lambda}}(\mu), $$
 where $(*)$ and $(**)$ follow respectively by Proposition \ref{prop:projective_shift} and  \eqref{eq:theta-projective_roots}. Hence, $\beta \in \mV_{\LL}(P_F)$ if and only if the one-dimensional $\K$-subspace $\K\cdot\mu$ of $\LL$ is contained in $\ker(F_{\alpha_\lambda})$.
 {However, by Lemma \ref{lem:normproperties}, the space $\xt(\ker(F_{\alpha_\lambda})\setminus\{0\})$ is a copy of a $(d_{\lambda}(F)-1)$-dimensional projective space over $\K$, and thus, so is $\mX_\lambda=\alpha_{\lambda}\cdot \xt(\ker(F_{\alpha_\lambda})\setminus\{0\})$.}
\end{proof}

\medskip

\subsection{Kernel and $\lambda$-Values}
Now we explore properties of the $\lambda$-values for a skew polynomial $F$, connecting them with the kernel of the endomorphisms associated to the skew polynomials $F_\alpha$'s.

\begin{proposition}\label{prop:equality_d_lambda}
  Let $F(X)\in \LL[X;\theta]$ be a nonzero skew polynomial and let $\alpha\in\LL^*$. Then, 
 $$\dim_\K(\ker F_\alpha)=d_{\Norm(\alpha)}(F).$$
\end{proposition}

\begin{proof}
 Let $\lambda\coloneqq \Norm(\alpha)$. By Corollary \ref{cor:kerf_skew}, $\ker(F_{\alpha})$ has dimension $\deg(\gcrd(F_{\alpha}(X),X^n-1))$, and we only need to show that this quantity is in turn equal to $d_{\lambda}(F)$. Clearly the ring homomorphism $G \longmapsto G_\alpha$  is degree-preserving. Hence, $\deg(\gcrd(F_{\alpha}(X),X^n-1))=\deg((\gcrd(F_{\alpha}(X),X^n-1))_{\alpha^{-1}})$, Moreover, by \eqref{eq:skew_homomoprhims_alpha}, we can deduce that $(\gcrd(F_{\alpha}(X),X^n-1))_{\alpha^{-1}}=\gcrd(F(X),(X^n-1)_{\alpha^{-1}})$. Finally, we obtain $(X^n-1)_{\alpha^{-1}}=\NN{\theta}{n}(\alpha^{-1})X^n-1=\lambda^{-1}X^n-1$ and thus $\deg(\gcrd(F_{\alpha}(X),X^n-1))=\deg(\gcrd(F(X),\lambda^{-1}X^n-1))=d_{\lambda}(F)$. 
\end{proof}

Here we introduce the following notation based on \cite{mcguire2019characterization}.  Let $F(X)=f_0+f_1X+\ldots+f_dX^d \in \LL[X;\theta]$ be a nonzero skew polynomial of degree $d$. Define the companion matrix $C_F\in \LL^{d\times d}$ as for classical commutative polynomials, and the matrix $A_F\in \LL^{d\times d}$ as
$$ A_F\coloneqq C_F  \theta(C_F)\cdot \ldots \cdot \theta^{n-1}(C_F)=C_FC_{X F X^{-1}}\cdot \ldots \cdot C_{X^{n-1} F  X^{-n+1}}.$$
Moreover, let us define the maps
$$ \begin{array}{rccl}
    \phi_F: &\LL[X;\theta]_{<d} &  \longrightarrow & \LL[X;\theta]_{<d} \\
    & G(X) & \longmapsto & X G(X) \mod\!\!_r\; F(X),
\end{array}$$

$$ \begin{array}{rccl}
    \psi_F: &\LL[X;\theta]_{<d} &  \longrightarrow & \LL[X;\theta]_{<d} \\
     &G(X) & \longmapsto & X^n G(X) \mod\!\!_r\: F(X),
\end{array}$$

where $\mod\!\!_r$ denotes the remainder obtained after the right division. The following result is due to McGuire and Sheekey in \cite{mcguire2019characterization}, where they prove it when $\LL$ is a finite field. However, since the proof adapts straightforwardly to our setting of a general field $\LL$, we leave it out.

\begin{proposition}
 Let $F\in \LL[X,\theta]$ be a nonzero skew polynomial of degree $d$ and let $\mB_d\coloneqq \{1, X, \ldots, X^{d-1}\}$. Then 
 \begin{enumerate}
     \item The map $\phi_F$ is a $\theta$-semilinear map of $\LL$-vector spaces whose associated matrix with respect to the basis $\mB_d$ is $A_F$.
     \item The map $\psi_F$ is an $\LL$-linear map whose associated matrix with respect to the basis $\mB_d$ is $C_F$.
 \end{enumerate}
\end{proposition}

With these new tools, we can deduce further connections between the kernel of the endomorphisms associated to the skew polynomials $F_\alpha$'s and the correspondent $\lambda$-values of $F$.

\begin{corollary}\label{cor:dim_kerf_alpha}
 Let $F(X)=f_0+\ldots+f_{d}X^{d}\in \LL[X;\theta]$ be a nonzero skew polynomial with $\deg (F)=d$ and let $\alpha\in\LL^*$. Then, we have
 $$\dim_\K(\ker F_\alpha)=\dim_\LL(\ker (A_F- \Norm(\alpha)I_d)).$$
\end{corollary}

\begin{proof}
 See \cite[Theorem 5]{mcguire2019characterization}
\end{proof}

We can finally derive the following important result, whose first part  was already shown in \cite[Proposition 1.3.7]{caruso2019residues}.

\begin{theorem}\label{thm:bound_degree_improved}
Let $F\in \LL[X;\theta]$ be a nonzero skew polynomial, let $\mathrm{A}\subseteq \LL^*$ be a set whose elements have pairwise distinct norms, and denote $\Lambda\coloneqq \Norm(\mathrm{A})$. Then
  the set $\{\alpha \in \mathrm{A} \mid  \ker(F_\alpha) \neq \{0\}\}$ is finite and 
  $$\sum_{\alpha \in \mathrm{A}} \dim \ker(F_\alpha)\leq \deg(F).$$
  Moreover, the following are equivalent:
  \begin{enumerate}
      \item $\sum_{\alpha \in \mathrm{A}} \dim \ker(F_\alpha)= \deg(F)$.
      \item $A_F$ is diagonalizable over $\LL$ and its eigenvalues are $\overline{\Lambda}\coloneqq \{\lambda \in \Norm(\mathrm{A})\mid   d_{\lambda}(F)>0\}$ where each $\lambda \in \overline{\Lambda}$ has multiplicity $d_{\lambda}(F)$.
      \item $F(X)$ right-divides the skew polynomial $H_{\Lambda}(X)\coloneqq \prod_{\lambda \in \Lambda}(X^n-\lambda)$.
  \end{enumerate}
\end{theorem}

\begin{proof}
 The inequality directly comes from Proposition \ref{prop:equality_d_lambda}, since we have
 \begin{align*}
     \sum_{\alpha \in \mathrm{A}} \dim \ker(F_\alpha)&=\sum_{\lambda \in \Norm(\mathrm{A})}d_{\lambda}(F) \\
     &=\sum_{\lambda\in\Norm(\mathrm{A})}\deg(\gcrd(F(X),X^n-\lambda)) \\
     &=\deg(\gcrd(F(X), H_{\Lambda}(X))\leq \deg(F),
 \end{align*}
 and equality occurs if and only if $F(X)$ right-divides $H_{\Lambda}(X)$.
 
 Furthermore, the size of the matrix $A_F$ is $d \times d$, where $d=\deg(F)$. Thus,  by Corollary \ref{cor:dim_kerf_alpha} we have 
 $$\sum_{\alpha \in \mathrm{A}} \dim \ker(F_\alpha)=\deg(F)$$
 if and only if 
 $$\sum_{\lambda \in \Norm(\mathrm{A})} \dim_\LL(\ker (A_F- \Norm(\alpha)I_d))=d,$$
 which in turn is equivalent to say that the matrix $A_F$ is diagonalizable over $\LL$ and its eigenvalues are $\overline{\Lambda}\coloneqq \{\lambda \in \Norm(\mathrm{A})\mid   d_{\lambda}(F)>0\}$, where each $\lambda \in \overline{\Lambda}$ has multiplicity $d_{\lambda}(F)$.
\end{proof}

\section{Skew Polynomial Framework for Sum-Rank Metric Codes} \label{sec:skewpoly_framework}
This section is devoted to explore a new point of view for sum-rank metric codes. We first show that we can define a isometries between the space $((\K^{n\times n})^\ell,\dd_{\MM})$. $(\LL^{\ell n},\dd_{\vv})$ and a suitable quotient of a skew polynomial ring. We then study sum-rank metric codes in this framework, with a particular focus on maximum sum-rank distance and linearized Reed-Solomon codes. 


For the whole section we fix the following notation. Let $\alpha_1,\ldots \alpha_\ell\in \LL^*$ be elements with pairwise distinct norms. Let  $\Lambda=\{\lambda_1,\ldots,\lambda_\ell\}\subseteq \K^*$ where $\lambda_i=\Norm(\alpha_i)$. Define 
$$H_\Lambda(X)=\prod_{i=1}^\ell (X^n-\lambda_i)\in \LL[X;\theta].$$
Observe that $H_\Lambda(X)$ belongs to the center of $\LL[X;\theta]$, that is $\K[X^n]$, since it is the product of $\ell$ central polynomials. Therefore, $H_\Lambda(X)$ generates a two-sided ideal in $\LL[X;\theta]$. 

\medskip

\subsection{Structure Theorem for Sum-Rank Metric Spaces}
In this section we study a new skew-polynomial framework in which we can naturally define the sum-rank metric.

The following result is a generalization of Theorem \ref{thm:ismorphism_skew_theta} and it provides a skew-polynomial description of the space $(\LL[\theta])^\ell$.

\begin{theorem}\label{thm:isomorphism_skew_sumrank}
 With the notation above, the map
 $$\begin{array}{rcl}\Phi_{\bs \alpha}: \LL[X;\theta] & \longrightarrow & (\LL[\theta])^\ell \\
 F(X) &\longmapsto &(\Phi(F_{\alpha_1}),\ldots, \Phi(F_{\alpha_\ell})) \end{array}$$
 is a surjective $\K$-algebra homomorphism, whose kernel is
 $(H_\Lambda(X)).$ Hence, it induces a $\K$-algebra isomorphism
 $$ \overline{\Phi}_{\bs \alpha}:\qspace \cong (\LL[\theta])^\ell.$$
 Moreover, both $\Phi_{\bs \alpha}$ and $\overline{\Phi}_{\bs \alpha}$ are also $\LL$-linear.
\end{theorem}

\begin{proof}
 One can see that the map $\Phi_{\bs\alpha}$ is obtained as the composition of the map
 $\rho:F \longmapsto (F_{\alpha_1,},\ldots,F_{\alpha_\ell})$ with the map defined as $\Phi$ acting entrywise. By Theorem \ref{thm:ismorphism_skew_theta},
 $\Phi$ is  a surjective $\K$-algebra homomorphism and an $\LL$-linear map. The $\LL$-linearity of $\rho$ is trivial. Moreover, the map $F\longmapsto F_\alpha$ is a ring homomorphism, due to \eqref{eq:skew_homomoprhims_alpha} and so it is $\rho$. Hence it is left to show that the map $\rho$ is surjective. For this purpose, we compute the kernel of $\Phi_{\bs\alpha}$. We have that $F\in\ker (\Phi_{\bs\alpha})$ if and only if $(X^n-1) \mid_{\rr}F_{\alpha_i}(X)$ for every $i\in\{1,\ldots,\ell\}$. Again, by \eqref{eq:skew_homomoprhims_alpha}, this is true if and only if $(X^n-1)_{\alpha_i^{-1}}\mid_{\rr}F(X)$. Since $(X^n-1)_{\alpha_i^{-1}}=\lambda_i^{-1}X^n-1$, we finally conclude that $F(X) \in \ker(\Phi_{\bs\alpha})$ if and only if $\lclm\{X^n-\lambda : \lambda \in \Lambda\}=H_\Lambda(X)$ right-divides $F(X)$. Thus, we have that $\mathrm{im}(\Phi_{\bs\alpha})\cong \qspace$. Since this is true also as $\LL$-vector spaces, and $\dim_{\LL}(\qspace)=\ell n=\dim_{\LL}((\LL[\theta])^\ell)$, we can conclude that $\Phi_{\bs\alpha}$ is in fact surjective, and $\overline{\Phi}_{\bs\alpha}$ is  an isomorphism.
 \end{proof}

\begin{remark}
 The isomorphism described in Theorem \ref{thm:isomorphism_skew_sumrank} clearly depends on the chosen elements $\alpha_1,\ldots,\alpha_\ell \in \LL^*$, and hence on the set $\Lambda$ of their norms. However, we can observe that   any set of $\ell$ elements $\alpha_1,\ldots,\alpha_\ell \in \LL^*$ whose norms are pairwise distinct induces an isomorphism.
\end{remark}

\begin{remark}
 It is easy to see that Theorem \ref{thm:ismorphism_skew_theta} is just a special case of Theorem \ref{thm:isomorphism_skew_sumrank} when $\ell=1$, and it is obtained by choosing $\alpha_1=1$.
\end{remark}

As a consequence of Theorem \ref{thm:isomorphism_skew_sumrank}, we can define the notion of sum-rank metric directly on the space $\qspace$.

\begin{definition}
 The \textbf{sum-rank weight} on the space $\qspace$ is the map 
 $$ \begin{array}{rccl} \wt_{\Lambda}:&\qspace & \longrightarrow & \N \\ &
 F(X) & \longmapsto & \ell n - \sum\limits_{i=1}^\ell d_{\lambda_i}(F)=\ell n-\deg \gcrd(F(X),H_\Lambda(X)).
 \end{array}$$
 Furthermore, the sum-rank weight induces the \textbf{sum-rank distance} on  $\qspace$,
 which is defined as
 $$\dd_{\Lambda}(F_1,F_2)\coloneqq \wt_{\Lambda}(F_1-F_2), \quad \mbox{ for any } F_1,F_2 \in \qspace.$$
\end{definition}

One can immediately ask whether this defintion of sum-rank weight and distance on $\qspace$ is consistent with the metric $\dd_{\srk}$  introduced for the space \red{$(\LL[\theta])^\ell$}. Not surprisingly, the answer is  that the map $\Phi_{\bs \alpha}$ is an isometry of metric spaces, as it shown in the following result.

\begin{corollary}\label{prop:isometry_skewpoly_skewalgebra} Let $\bs\alpha=(\alpha_1,\ldots,\alpha_n)\in\LL^\ell$ be any vector such that $\Norm(\alpha_i)=\lambda_i$ for each $i=1,\ldots,\ell$. For every $F, G \in \qspace$, we have 
 $$\dd_{\Lambda}(F,G)=\dd_{\srk}(\overline{\Phi}_{\bs \alpha}(F),\overline{\Phi}_{\bs \alpha}(G)).$$
 In particular, the map $\overline{\Phi}_{\bs \alpha}:(\qspace, \dd_{\Lambda})\rightarrow((\LL[\theta])^{\ell},\dd_{\srk}))$ is an isometry of metric spaces.
\end{corollary}

\begin{proof}
 It directly comes from the definition of the two metrics together with Proposition \ref{prop:equality_d_lambda}.
\end{proof}

\medskip

\subsection{The $\Lambda$-Dickson Matrix}

We now show that, as in the case of rank-metric codes and their representations as linearized polynomials or skew group algebra elements, it is possible to relate the sum-rank weight of a skew polynomial to the rank of a particular matrix of Dickson-type; see e.g.  \cite{menichetti1986roots,wu2013linearized} for finite fields and \cite{augot2020rank} for general Galois extensions.

For this purpose, for a given   $F\in\qspace$, define  the $\LL$-linear map
$$\mu: \qspace \longrightarrow \End_{\LL}(\qspace),$$
given by
$$ \begin{array}{rccl}\mu(F):&
  \qspace & \longrightarrow &\qspace, \\ & G(X) & \longmapsto & G(X)  F(X). \end{array}$$

\begin{definition}
 Let $F\in\qspace$. The matrix representation of the $\LL$-linear map $\mu(F)$ with respect to the monomial basis $\{1,X,\ldots,X^{\ell n-1}\}$ is called the \textbf{$\Lambda$-Dickson matrix} associated to $F$, and it is denoted by $D_\Lambda(F)\in\LL^{\ell n \times \ell n}$.
\end{definition}

 Dickson matrices associated to linearized polynomials over finite fields are well-known objects \cite{menichetti1986roots}.   Such matrices have been recently generalized to any Galois extension of fields in \cite{augot2020rank}, where they were called \emph{$G$-Dickson matrices}. We kept a similar name because when $\Lambda$ is a finite subgroup, the $\Lambda$-Dickson matrix resembles these objects, as the following example shows.

\begin{example}
    Let $\Lambda$ be a finite -- and hence cyclic -- subgroup of $\K^*$ of order $\ell$ and let us consider  $F(X) = f_0 +\ldots+f_{\ell n-1}X^{\ell n-1} \in \qspace$. Since $\Lambda$ is cyclic, we have $H_{\Lambda}(X)=X^{\ell n}-1$. In this case, it is not difficult to see that  the  $\Lambda$-Dickson matrix
  associated to $F(X)$ is $D_\Lambda(F)=(d_{i,j})\in \LL^{\ell n \times \ell n}$, where 
 \[
 d_{i,j} = \theta^{j-1}(f_{(j-i \mod \ell n)}), \quad\quad \forall\ i,j \in
 \{1, \dots, \ell n\}.
 \]
\end{example}

\begin{definition} Let $F\in \qspace$, The (left) \textbf{annihilator} of the skew polynomial $F$ is 
$$\Ann (F)\coloneqq  \{G \in \qspace  \mid G(X)F(X)=0\}.$$
\end{definition}

It is easy to see that the left annihilator of a skew polynomial in $\qspace$ is a left ideal, and it is clearly principal. Moreover, one has an explicit formula for the generator of such a left ideal.

\begin{proposition}[\textnormal{\cite[Theorem 14]{ore1933theory}}]\label{prop:annihilator} Let $F\in \qspace$. Then
  $$\Ann(F)=(\qspace)(A(X)),$$
  where $A(X)\coloneqq \lclm(H_\Lambda(X),F(X))F(X)^{-1}$.\footnote{Since $\lclm(H_\Lambda(X),F(X))$ is right divisible by $F(X)$, then  $\lclm(H_\Lambda(X),F(X))=B(X)F(X)$  for some skew polynomial $B(X)$. Thus, one defines $\lclm(H_\Lambda(X),F(X))F(X)^{-1}\coloneqq B(X)$.}
\end{proposition}

We are now ready to connect the sum-rank weight of a skew polynomial with the rank of its associated Dickson matrix.

\begin{theorem}
 For every $F(X)\in\qspace$ it holds that
 $$\wt_\Lambda(F)=\ell n-\dim_{\LL}(\Ann(F))=\rk(D_\Lambda(F)).$$
\end{theorem}

\begin{proof}
 By Proposition \ref{prop:annihilator}, we have that $\Ann(F)$ is the left ideal generated by $A(X)=\lclm(H_\Lambda(X),F(X))F(X)^{-1}$. Moreover, it not difficult to see that  $A(X) \mid_{\rr} H_\Lambda(X)$ (see e.g. \cite[Theorem  11]{ore1933theory}). Therefore, 
 \begin{align*}\dim_{\LL}(\Ann(F))&=\deg(H_\Lambda)-\deg(A)=\deg(H_\Lambda)-\deg(\lclm(H_\Lambda,F))+\deg(F)\\
 &=\deg(H_\Lambda)-\deg(H_\Lambda)-\deg(F)+\deg(\gcrd(H_\Lambda,F))+\deg(F)\\ &=\deg(\gcrd(H_\Lambda,F)),
 \end{align*}
 which shows the first equality. For the second equality, it is enough to observe that $\Ann(F)=\ker(\mu(F))$, and then we conclude.
\end{proof}

\begin{example}\label{exa:f5}
Let $\K=\F_5$, $n=3$, $\ell=4$ and let $\theta$ be the Frobenius automorphism defined by $\theta(\alpha)=\alpha^5$. Since $\gcd(n,q-1)=1$, we have that $\F_5^*$ is a set of representatives for the norm function. Hence, we can take $\alpha_i=i$  and $\lambda_i=\Norm(\alpha_i)=i^3$ for $i\in\{1,2,3,4\}$. Thus, we have $\Lambda=\F_5^*$ and we can compute the polynomial 
 $$H_\Lambda(X)=\prod_{i=1}^4(X^3-i)=X^{12}-1.$$
 Now let us take $F(X)=X^4+2X^3+3X^2+3X+1 \in \F_{5^3}[X;\theta]/(X^{12}-1)$. With straightforward computations one can see that $F(X)$ right-divides $X^{12}-1$ and we have $d_1(F)=1$, $d_2(F)=0$, $d_3(F)=2$ and $d_4(F)=1$. Therefore,
 $$\wt_\Lambda(F)=\ell n-\deg(\gcrd(F,H_\Lambda))=12-4=8.$$
 The $\ell$ components obtained through the map $\overline{\Phi}_{\bs\alpha}$ are given by
 $$\begin{array}{rclrcl}
    \overline{\Phi}(F_1)& \!\!\!\!=\!\!\!\!& 3\theta^2 + 4\theta + 3 \mathrm{id}, \;\;\; &\;\;\;\overline{\Phi}(F_2)&\!\!\!\!=\!\!\!\!&2\theta^2+2\theta+2  \mathrm{id} \\
    \overline{\Phi}(F_3)&\!\!\!\!=\!\!\!\!& 2\theta^2, & \overline{\Phi}(F_4)&\!\!\!\!=\!\!\!\!& 3\theta^2 + 3\theta + 4 \mathrm{id}
 \end{array}$$
 and their rank-weights are 
 \begin{align*} 
 \rk(\overline{\Phi}(F_1))&=n-d_{\Norm(1)}(F)=3-d_1(F)=2, \\
  \rk(\overline{\Phi}(F_2))&=n-d_{\Norm(2)}(F)=3-d_3(F)=1, \\
   \rk(\overline{\Phi}(F_3))&=n-d_{\Norm(3)}(F)=3-d_2(F)=3, \\
    \rk(\overline{\Phi}(F_4))&=n-d_{\Norm(4)}(F)=3-d_4(F)=2.
 \end{align*}
 Furthermore, we can compute $\Ann(F)$, which by Proposition \ref{prop:annihilator} is the  left ideal generated by the polynomial
 $$(\lclm(H_\Lambda(X),F(X))F(X)^{-1}= x^8 + 3x^7 + x^6 + x^5 + x^3 + 4x^2 + 3x + 4.$$
 Finally, the $\Lambda$-Dickson matrix associated to $F$ is
$$D_\Lambda(F)=\left(\begin{array}{cccccccccccc}
1 & 0 & 0 & 0 & 0 & 0 & 0 & 0 & 1 & 2 & 3 & 3 \\
3 & 1 & 0 & 0 & 0 & 0 & 0 & 0 & 0 & 1 & 2 & 3 \\
3 & 3 & 1 & 0 & 0 & 0 & 0 & 0 & 0 & 0 & 1 & 2 \\
2 & 3 & 3 & 1 & 0 & 0 & 0 & 0 & 0 & 0 & 0 & 1 \\
1 & 2 & 3 & 3 & 1 & 0 & 0 & 0 & 0 & 0 & 0 & 0 \\
0 & 1 & 2 & 3 & 3 & 1 & 0 & 0 & 0 & 0 & 0 & 0 \\
0 & 0 & 1 & 2 & 3 & 3 & 1 & 0 & 0 & 0 & 0 & 0 \\
0 & 0 & 0 & 1 & 2 & 3 & 3 & 1 & 0 & 0 & 0 & 0 \\
0 & 0 & 0 & 0 & 1 & 2 & 3 & 3 & 1 & 0 & 0 & 0 \\
0 & 0 & 0 & 0 & 0 & 1 & 2 & 3 & 3 & 1 & 0 & 0 \\
0 & 0 & 0 & 0 & 0 & 0 & 1 & 2 & 3 & 3 & 1 & 0 \\
0 & 0 & 0 & 0 & 0 & 0 & 0 & 1 & 2 & 3 & 3 & 1 \\
\end{array}\right)\in (\F_{5^3})^{12 \times 12},$$
which has rank $8$.
\end{example}

\begin{example}\label{exa:f5new}
We fix the same setting as the one in Example \ref{exa:f5} and choose the primitive element $\gamma$ to be a root of $y^3+3y+3$. If we  select the skew polynomial $G(X)=X^4+\gamma^{55}X^3+\gamma^{29}X^2+\gamma^{63} X+1$ and compute the $\lambda$-values, we obtain 
$$ d_1(G)=d_3(G)=1, \qquad d_2(G)=d_4(G)=0.$$
Thus, we have $\wt_\Lambda(G)=12-2=10$. Its $\Lambda$-Dickson matrix is given by
$$D_\Lambda(G)=\small{\left(\begin{array}{cccccccccccc}
1 & 0 & 0 & 0 & 0 & 0 & 0 & 0 & 1 & \gamma^{55} & \gamma^{21} & \gamma^{87} \\
\gamma^{63} & 1 & 0 & 0 & 0 & 0 & 0 & 0 & 0 & 1 & \gamma^{27} & \gamma^{105} \\
\gamma^{29} & \gamma^{67} & 1 & 0 & 0 & 0 & 0 & 0 & 0 & 0 & 1 & \gamma^{11} \\
\gamma^{55} & \gamma^{21} & \gamma^{87} & 1 & 0 & 0 & 0 & 0 & 0 & 0 & 0 & 1 \\
1 & \gamma^{27} & \gamma^{105} & \gamma^{63} & 1 & 0 & 0 & 0 & 0 & 0 & 0 & 0 \\
0 & 1 & \gamma^{11} & \gamma^{29} & \gamma^{67} & 1 & 0 & 0 & 0 & 0 & 0 & 0 \\
0 & 0 & 1 & \gamma^{55} & \gamma^{21} & \gamma^{87} & 1 & 0 & 0 & 0 & 0 & 0 \\
0 & 0 & 0 & 1 & \gamma^{27} & \gamma^{105} & \gamma^{63} & 1 & 0 & 0 & 0 & 0 \\
0 & 0 & 0 & 0 & 1 & \gamma^{11} & \gamma^{29} & \gamma^{67} & 1 & 0 & 0 & 0 \\
0 & 0 & 0 & 0 & 0 & 1 & \gamma^{55} & \gamma^{21} & \gamma^{87} & 1 & 0 & 0 \\
0 & 0 & 0 & 0 & 0 & 0 & 1 & \gamma^{27} & \gamma^{105} & \gamma^{63} & 1 & 0 \\
0 & 0 & 0 & 0 & 0 & 0 & 0 & 1 & \gamma^{11} & \gamma^{29} & \gamma^{67} & 1 \\
\end{array}\right)}\in (\F_{5^3})^{12 \times 12},$$
\end{example} 

We can now introduce the notion of codes endowed with the sum-rank metric in this skew polynomial framework.

\begin{definition}
 A  \textbf{sum-rank metric code} $\C$ is a subset of $\qspace$ endowed with the sum-rank distance $\dd_\Lambda$. The \textbf{minimum sum-rank distance} of $\C$ is the integer
 $$ \dd_{\Lambda}(\C)\coloneqq \min\{\dd_\Lambda(F,G) \mid F,G \in \C, F\neq G\}.$$
 
 Moreover,  for any  subfield $\F\subseteq \LL$, a sum-rank metric code will be called $\F$-linear, if it is an $\F$-linear subspace of $\qspace$. 
 \end{definition}

\medskip

\subsection{Connections with Classical Frameworks}\label{sec:comparison_frameworks}

Here we briefly recap how to recover the original frameworks of sum-rank metric codes, seen either as spaces of vectors of length $\ell n$ over $\LL$, or as subspaces of the direct sum of $\ell$ matrix spaces over $\K$.

First, let $\mE$ be a $\K$-basis of $\LL$. Write every element of
$\LL$ in coordinates with respect to $\mE$, resulting in a column
vector in $\K^n$.  In the same way, we can transform a vector
$v \in \LL^n$ to a matrix in $\K^{n \times n}$, which we denote by
$\Ext_\mE(v)$.
Let $\mB_1,\ldots,\mB_\ell$ be $\K$-bases of $\LL$ and write $\bs\mB\coloneqq (\mB_1,\ldots,\mB_\ell)$. If $v=(v^{(1)} \mid \cdots \mid v^{(\ell)})\in \LL^{\ell n}$, we define
$$ \Ext_{\mB}(v)\coloneqq (\Ext_{\mB_1}(v^{(1)}) , \ldots, \Ext_{\mB_\ell}(v^{(\ell)}))\in(\mat)^\ell.$$
Given a vector sum-rank metric code $\C\subseteq \LL^{\ell n}$, we define
$$\qquad \qquad\qquad  \qquad \qquad \Ext_{\bs\mB}(\C)\coloneqq \{\Ext_{\bs\mB}(c) \mid c\in\C \}. \qquad \qquad \qquad\qquad(\red{\mathrm{vec}})\rightarrow (\red{\mathrm{mat}})$$
From the definition of the metrics, it is immediate to verify that the map $\Ext_{\mB}:(\LL^{\ell n},\dd_{\vv}) \rightarrow ((\mat)^{\ell},\dd_{\MM})$ is an isometry, and hence one can equivalently study codes in one or in the other setting.  

\medskip

Now, for a given $\K$-basis $\mE=(e_1,\ldots,e_n)$ of $\LL$, define the map
$$\begin{array}{rcl}
\ev_{\mE}: \LL[\theta] & \longrightarrow & \LL^n \\
f & \longmapsto & (f(e_1),\ldots,f(e_n)). 
\end{array}$$
Denote by 
$\ev_{\bs\mB}: (\LL[\theta])^{\ell}\rightarrow(\LL^{\ell n})$, the direct sum of the maps $\ev_{\mB_i}$.
Given a sum-rank metric code $\C\subseteq \qspace$, define
$$\qquad \qquad\qquad \qquad \qquad \qquad  \C(\bs\alpha,\bs\mB)\coloneqq \ev_{\bs\mB}(\overline{\Phi}_{\bs\alpha}(\C)). \qquad \qquad\qquad \qquad\qquad(\red{\mathrm{poly}})\rightarrow (\red{\mathrm{vec}}) $$
Also in this case it is straightforward to verify that
$$\ev_{\bs\mB}\circ \overline{\Phi}_{\bs\alpha}:(\qspace,\dd_{\vv}) \rightarrow (\LL^{\ell n},\dd_{\vv})$$
is an isometry of metric spaces. This is due to the fact that each map $\ev_{\mB_i}$ is an isometry (see e.g. \cite{augot2020rank}) -- and hence also the direct sum of them -- and to the fact that $\overline{\Phi}_{\bs\alpha}$ is an isometry by Corollary \ref{prop:isometry_skewpoly_skewalgebra}.

\medskip

Finally, one can combine the isometries described above. Suppose that we have $\bs\mB=(\mB_1,\ldots,\mB_\ell), \bs\mE=(\mE_1,\ldots,\mE_\ell)$, where the $\mB_i$'s and the $\mE_i$'s are $\K$-bases of $\LL$. For a given sum-rank metric code $\C\subseteq \qspace$ we define
$$\qquad \qquad\qquad\qquad\qquad \Mat_{\bs\alpha,\bs\mB,\bs\mE}(\C)\coloneqq \Ext_{\bs\mE}(\C(\bs\alpha,\bs\mB)). \qquad\qquad\qquad \qquad(\red{\mathrm{poly}})\rightarrow (\red{\mathrm{mat}})$$
Since the map $\Mat_{\bs\alpha,\bs\mB,\bs\mE}:(\qspace,\dd_{\Lambda}) \rightarrow ((\mat)^\ell,\dd_{\MM})$ is the composition of two isometries, is itself an isometry of metric spaces.

\medskip
The procedures described above show how to get a code representation from another. Moreover, all the procedures are  reversible, assuming that there exists a cyclic Galois extension of $\K$ of degree $n$. More specifically,  from a skew-polynomial representation of sum-rank metric codes we can always derive a vector representation, and from a vector representation we can always get a matrix representation. On the other hand, from a matrix sum-rank metric code over $\K$ we can derive a vector representation only when $\K$ admits a degree $n$ extension field $\LL$, and the vector representation can be converted into a skew-polynomial representation as $\qspace$, only when the extension $\LL/\K$ is Galois and cyclic and $\ell \leq |\K^*|$.

\medskip
\begin{center}
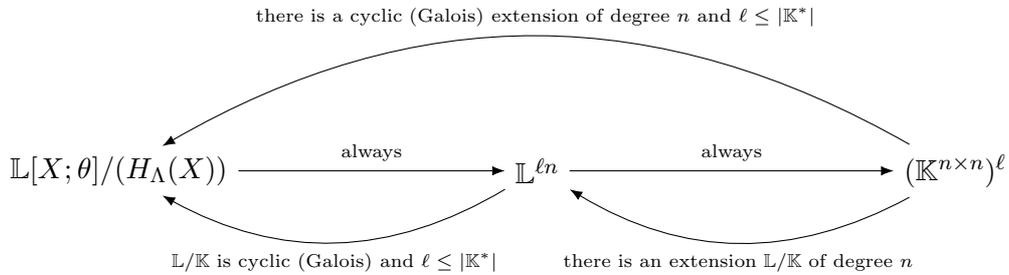
\begin{figure}[!h]
\begin{tikzpicture}[node distance = 5.5cm, auto]
      \node (K) {$\qspace$};
      \node (F) [right of=K] {$\LL^{\ell n}$};
      \node (E) [right  of=F] {$(\mat)^\ell$};
      \draw[-Latex] (K) to node {\tiny{always}} (F);
      \draw[-Latex] (F) to node {\tiny{always}} (E);
      \draw[-Latex, bend right] (E) to node [swap] {\tiny{there is a cyclic (Galois) extension of degree $n$  and $\ell \leq |\K^*|$}} (K);
      \draw[-Latex, bend left] (E) to node  {\tiny{there is an extension $\LL/\K$ of degree $n$}} (F);
      \draw[-Latex, bend left] (F) to node  {\tiny{$\LL/\K$ is cyclic (Galois) and $\ell \leq |\K^*|$}} (K);
      \end{tikzpicture}
      \caption{The diagram above shows under which conditions it is possible to transform a framework into another one.}\label{fig:1}
      \end{figure}
\end{center}
However, over finite fields -- that is where most applications come from -- the three representations are equivalent whenever $\ell \leq q-1$, due to the fact that for every positive integer $n$ there always exists a degree $n$-extension of a finite field $\Fq$ and such an extension field is also cyclic. When we choose  $\K=\mathbb Q$ the three representations are completely equivalent, since it is well-known that for every integer $n$ there exists a cyclic Galois extension of degree $n$ over $\mathbb Q$. On the other hand, when $\K=\mathbb R$, the three representations are equivalent only when $n=2$, while for $n\geq 3$ there is no vector representation, and hence neither a skew-polynomial representation.

\begin{example}\label{exa:f5bis} Let us consider the same setting and notation of Examples \ref{exa:f5} and \ref{exa:f5new}: let $\K=\F_5$, $n=3$, $\ell=4$ and let $\theta$ be the $5$-Frobenius automorphism. Moreover, let us take $\alpha_i=i$  and $\lambda_i=\Norm(\alpha_i)=i^3$ for $i\in\{1,2,3,4\}$. Thus, we fox $\Lambda=\F_5^*$ and 
 $$H_\Lambda(X)=\prod_{i=1}^4(X^3-i)=X^{12}-1.$$
 Now consider the primitive element $\gamma$ of $\F_{5^3}$ that is a root of the polynomial $y^3+3y+3$ and let us choose  $F(X)=X^4+2X^3+3X^2+3X+1 \in \F_{5^3}[X;\theta]/(X^{12}-1)$. At this point, define the code
 $$ \C\coloneqq \langle 1, F(X)\rangle_{\F_{5^3}}.$$
 As $\F_5$-bases of $\F_{5^3}$ we take $\bs\mB\coloneqq (\mathcal B, \mathcal B, \mathcal B,\mathcal B)$, where $\mathcal B\coloneqq (1,\gamma,\gamma^2)$. We have already computed $\overline{\Phi}_{\bs\alpha}(F)$ in Example \ref{exa:f5}, and its components are 
 $$\begin{array}{rclrcl}
    \overline{\Phi}(F_1)& \!\!\!\!=\!\!\!\!& 3\theta^2 + 4\theta + 3 \mathrm{id}, \;\;\; &\;\;\;\overline{\Phi}(F_2)&\!\!\!\!=\!\!\!\!&2\theta^2+2\theta+2  \mathrm{id}, \\
    \overline{\Phi}(F_3)&\!\!\!\!=\!\!\!\!& 2\theta^2, & \overline{\Phi}(F_4)&\!\!\!\!=\!\!\!\!& 3\theta^2 + 3\theta + 4 \mathrm{id}.
 \end{array}$$
 Moreover, we clearly have $\overline{\Phi}_{\bs\alpha}(1)=(\mathrm{id},\mathrm{id},\mathrm{id},\mathrm{id})$. Thus, 
 \begin{align*} \C(\bs\alpha,\bs\mB)&=\langle \evB(\overline{\Phi}_{\bs\alpha}(1)),\evB(\overline{\Phi}_{\bs\alpha}(F)) \rangle_{\F_{5^3}}  \\
 &= \mathrm{Rowsp}\left(\begin{array}{ccc|ccc|ccc|ccc} 1 & \gamma & \gamma^2 & 1 & \gamma & \gamma^2 & 1 & \gamma & \gamma^2 & 1 & \gamma & \gamma^2 \\
 0 & \gamma^5 & \gamma^{32} & 1 & 0 & \gamma^{93} & \gamma^{31} & \gamma^{56} & \gamma^{81} & 0 & \gamma & \gamma^{56} \\ 
 \end{array}\right).
 \end{align*}
 
 Now, let us fix also $\bs\mE\coloneqq \bs\mB=(\mB,\mB,\mB,\mB)$ and let us denote by $u$ and $v$ the first and respectively the second row of the matrix above. With easy computations one gets 
 $$ \Ext_{\bs\mB}(u)=\left(\begin{pmatrix} 1 & 0 & 0 \\ 0 & 1 & 0 \\ 0 & 0 & 1\end{pmatrix},\begin{pmatrix} 1 & 0 & 0 \\ 0 & 1 & 0 \\ 0 & 0 & 1\end{pmatrix},\begin{pmatrix}1 & 0 & 0 \\ 0 & 1 & 0 \\ 0 & 0 & 1 \end{pmatrix},\begin{pmatrix} 1 & 0 & 0 \\ 0 & 1 & 0 \\ 0 & 0 & 1\end{pmatrix}\right), $$
  $$ \Ext_{\bs\mB}(v)=\left(\begin{pmatrix} 0 & 4 & 0 \\ 0 & 4 & 2 \\ 0 & 2 & 0\end{pmatrix},\begin{pmatrix}1 & 0 & 3 \\ 0 & 0 & 0 \\ 0 & 0 & 0 \end{pmatrix},\begin{pmatrix} 2 & 2 & 2 \\ 0 & 0 & 1 \\ 0 & 1 & 3\end{pmatrix},\begin{pmatrix} 0 & 0 & 2 \\ 0 & 1 & 0 \\ 0 & 0 & 1\end{pmatrix}\right). $$
  If  we define $M$ to be the companion matrix of the polynomial $y^3+3y+3$, that is 
  $$ M=\begin{pmatrix} 0 & 0 & 2 \\ 1 & 0 & 2 \\ 0 & 1 & 0
  \end{pmatrix},$$
  we then obtain that
  $$ \Mat_{\bs\alpha,\bs\mB,\bs\mB}(\C)=\langle \{M^i \cdot  \Ext_{\bs\mB}(u) : 0 \leq i \leq 2\} \cup \{ M^i\cdot\Ext_{\bs\mB}(v) : 0 \leq i \leq 2\}\rangle_{\F_5}.$$
\end{example}

\medskip

\subsection{The Adjoint Operator}
In the theory of $\theta$-polynomials, the adjoint  defined in \eqref{eq:adjoint} plays an important role in $\LL[\theta]$: it provides in this setting the analogue of the transpose in $\K^{n \times n}$. If we want to understand the analogue for the ambient space $(\K^{n\times n})^{\ell}$, then we may transpose   any $r$ matrices in an $\ell$-uple of matrices. This can be formalized via the action of the group $(\ZZ{2})^\ell$ on $(\K^{n\times n})^\ell$, given by
$$ (v_1,\ldots,v_\ell)\cdot (A_1,\ldots, A_\ell) =(B_1,\ldots,B_\ell),$$
where
$$B_i=\begin{cases} A_i & \mbox{ if } v_i=0, \\
A_i^\top & \mbox{ if } v_i=1. 
\end{cases}$$

Here we give the counterpart to this map in our setting. First, for any vector $v=(v_1,\ldots,v_\ell) \in (\ZZ{2})^\ell$ define the map $\Upsilon_v: (\LL[\theta])^\ell\longrightarrow (\LL[\theta])^\ell$ as
$$ \Upsilon_v(f_1,\ldots, f_\ell) =(g_1,\ldots,g_\ell),$$
where
$$g_i=\begin{cases} f_i & \mbox{ if } v_i=0, \\
f_i^\top & \mbox{ if } v_i=1. 
\end{cases}$$

\begin{remark}\label{rem:well_defined_phibeta}
 Observe that for any set of elements $\mathrm{A}'=\{\beta_1,\ldots, \beta_\ell\}$ such that $\Norm(\mathrm{A}')=\Norm(\mathrm{A})=\Lambda$, by Theorem \ref{thm:isomorphism_skew_sumrank} the map $\overline{\Phi}_{\bs\beta}:(\qspace,\dd_\Lambda)\longrightarrow ((\LL[\theta))^\ell,\dd_{\srk})$ is a well-defined isometry, where $\bs\beta:=(\beta_1,\ldots,\beta_\ell)$.
\end{remark}

\begin{definition}
 Let $\bs\alpha=(\alpha_1,\ldots,\alpha_\ell)$, $\bs\beta=(\beta_1,\ldots,\beta_\ell) \in \LL^{\ell}$ be two vectors such that $\Norm(\{\alpha_1,\ldots,\alpha_{\ell}\})=\Norm(\{\beta_1,\ldots,\beta_{\ell}\})=\Lambda$.  For any $v\in (\ZZ{2})^\ell$, and any $F(X)\in\qspace$, the \textbf{$v$-adjoint} of $F(X)$ with respect to $\bs\alpha$ and $\bs\beta$ is the skew polynomial
 $$ F^{v}_{[\bs\alpha,\bs\beta]}(X):=(\overline{\Phi}_{\bs\beta}^{-1} \circ \Upsilon_v \circ \overline{\Phi}_{\bs\alpha})(F(x)).$$
\end{definition}

Clearly, one has
$$ \wt_{\Lambda}(F^v_{[\bs\alpha,\bs\beta]})=\wt_{\Lambda}(F), \qquad \mbox{ for every } F \in \qspace \mbox{ and every } v \in (\ZZ{2})^\ell.$$
Note that it seems  difficult to express the $v$-adjoint of an element $F(X)$ in a simple formula using the coefficients of $F(X)$. However, there is a special case in which we can derive a nice expression.

 For a $\bs\alpha=(\alpha_1,\ldots,\alpha_\ell)$, we will write $\bs\alpha^{-1}$ to indicate the vector whose entries are the inverses of the entries of $\bs\alpha$, that is $\bs\alpha^{-1}\coloneqq (\alpha_1^{-1},\ldots, \alpha_{\ell}^{-1})$. 
 In particular, if we assume that $\Lambda$ is a subgroup, we have that $\Norm(\{\alpha_1^{-1},\ldots,\alpha_\ell^{-1}\})=\Norm(\mathrm{A})=\Lambda$, and 
$$\overline{\Phi}_{\bs\alpha^{-1}}:\qspace \longrightarrow (\LL[\theta])^\ell $$
is well defined by Remark \ref{rem:well_defined_phibeta}.

\begin{theorem}\label{thm:adjoint}
 Let $\mathbf{1}\in(\ZZ{2})^\ell$ denote the all-ones vector, and let $F(X)=f_0+f_1X+\ldots +f_{\ell n-1} X^{\ell n-1}\in\qspace$. Assume that $\Lambda$ is a cyclic group of order $\ell$. Then
 $$F^{\mathbf{1}}_{[\bs\alpha,\bs\alpha^{-1}]}(X)=\sum_{i=1}^{\ell n}\theta^{n-i}(f_i)X^{\ell n-i}.\footnote{Here we write the indices of the coefficients  modulo $\ell n$. In particular, in the formula above, we define $f_{\ell n}:=f_0$.} $$
\end{theorem}

\begin{proof}
 Notice that it is enough to show that $(\overline{\Phi} (\mu(\alpha X)^u)^\top=\overline{\Phi}(\theta^{n-u}(\mu)(\alpha^{-1}X)^{\ell n-u})$, for every $u\in\{0,\ldots,\ell n-1\}$, $\alpha,\mu \in \LL^*$ with $\Norm(\alpha)\in\Lambda$, since then we can conclude by $\K$-linearity. Let $\lambda\coloneqq\Norm(\alpha)$, and let us write $u=i+tn$, with $0\leq i \leq n-1$, and $0\leq t\leq \ell-1$. Then 
 \begin{align*}
     (\overline{\Phi} (\mu(\alpha X)^u)^\top&=(\mu\NN{\theta}{i}(\alpha)\lambda^t\theta^i)^\top  \\
     &=\theta^{n-i}(\mu)\NN{\theta}{i}(\theta^{n-i}(\alpha))\theta^{n-i} \\
     &= \theta^{n-i}(\mu)\lambda^{t+1}\NN{\theta}{n-i}(\alpha^{-1})\theta^{n-i} \\
     &= \theta^{n-i}(\mu)(\lambda^{-1})^{\ell-t-1}\NN{\theta}{n-i}(\alpha^{-1})\theta^{n-i}.\\
     &=\overline{\Phi}(\theta^{n-u}(\mu)(\alpha^{-1}X)^{(\ell-t-1)n+n-i}) \\
     &=\overline{\Phi}(\theta^{n-u}(\mu)(\alpha^{-1}X)^{\ell n-u}).
 \end{align*}
\end{proof}

In particular, when $\Lambda$ is a cyclic group, we deduce by Proposition \ref{thm:adjoint} that the $\mathbf{1}$-adjoint $F^{\mathbf{1}}_{[\bs\alpha,\bs\alpha^{-1}]}(x)$ does not depend on $\bs\alpha$. Thus, from now on, we will write
$$ F^\top(X):=F^{\mathbf{1}}_{[\bs\alpha,\bs\alpha^{-1}]}(X).$$

Using the formula provided in Theorem \ref{thm:adjoint}, we can immediately deduce the following result on the $\Lambda$-Dickson matrix of the $\mathbf{1}$-adjoint of a skew polynomial in $\qspace$.
\begin{corollary}
 Let $\Lambda$ be a cyclic group, and let $F(X)\in\qspace$. Then
 $$ D_{\Lambda}(F^\top)=D_{\Lambda}(F)^\top.$$
\end{corollary}

We conclude this section by providing the definition of adjoint of a code.

\begin{definition}
 Let $\bs\alpha=(\alpha_1,\ldots,\alpha_\ell)$, $\bs\beta=(\beta_1,\ldots,\beta_\ell) \in \LL^{\ell}$ be two vectors such that $\Norm(\{\alpha_1,\ldots,\alpha_{\ell}\})=\Norm(\{\beta_1,\ldots,\beta_{\ell}\})=\Lambda$, and let $v\in (\ZZ{2})^\ell$. The \textbf{$v$-adjoint code} of a $\K$-linear code $\C\subseteq \qspace$ is
 $$ \C_{[\bs\alpha,\bs\beta]}^v:=\left\{F_{[\bs\alpha,\bs\beta]}^v(X) \colon F(X)\in \C\right\}.$$
 Furthermore, when $v=\mathbf{1}$, $\Lambda$ is a cyclic group and $\bs\beta=\bs\alpha^{-1}$, we will simply refer to it as the \textbf{adjoint code} and we will denote it by $\C^{\top}$.
\end{definition}

\medskip

\subsection{Maximum Sum-Rank Distance and Linearized Reed-Solomon Codes} 
We now focus on the special class of sum-rank metric codes which have optimal parameters. 
The following is a well-known result, already proved in \cite{martinez2018skew} for vector codes.

\begin{theorem}\label{thm:singleton}
 Let $\F$ be a subfield of $\LL$, and let $\C\subseteq \qspace$ be an $\F$-linear sum-rank metric code. Then 
 $$\dim_{\F}(\C)\leq [\LL:\F](\ell n-\dd_{\Lambda}(\C)+1).$$
\end{theorem}

\begin{proof}
 We first identify the codewords as vectors in $\LL^{\ell n}$ via $\ev_{\bs\mB} \circ \overline{\Phi}_{\bs\alpha}$. In other words, we consider the vector code $\C(\bs\alpha,\bs\mB)$, for some suitable vector $\bs\alpha$ whose entries have norm equal to $\lambda_i$'s and any vector $\bs\mB$ formed by $\ell$ bases of $\LL$ over $\K$. We will prove that the $\F$-dimension of $\C(\bs\alpha,\bs\mB)$ is at most $[\LL:\F](\ell n-\dd_{\Lambda}(\C)+1)$. Let us consider the \emph{Hamming metric} on $\LL^{\ell n}$, defined via the Hamming weight $\wt_{\HH}(v)\coloneqq|\{i: v_i \neq 0\}|$. Let $\dd_{\HH}$ denote the minimum Hamming distance of  $\C(\bs\alpha,\bs\mB)$. It is immediate to see that 
 $$\dd_{\HH}\geq \dd_{\vv}(\C(\bs\alpha,\bs\mB))=\dd_{\Lambda}(\C).$$
 If we now puncture the code $\C(\bs\alpha,\bs\mB)$ on $\dd_{\HH}-1$ coordinates, then by definition of Hamming distance, the obtained code has the same $\F$-dimension. However, this code is contained in the space $\LL^{\ell n-\dd_{\HH}+1}$, which has $\F$-dimension equal to $[\LL:\F](\ell n-\dd_{\HH}+1)$. Thus, we deduce $$\dim_{\F}(\C)\leq [\LL:\F](\ell n-\dd_{\HH}(\C)+1)\leq [\LL:\F](\ell n-\dd_{\Lambda}(\C)+1).$$
\end{proof}

\begin{definition}
 An $\F$-linear code $\C\subseteq \qspace$ is called \textbf{maximum sum-rank distance} (or \textbf{MSRD} in short) if its parameters meet the Singleton-like bound of Theorem \ref{thm:singleton} with equality. 
\end{definition}

 The study of MSRD codes has started in \cite{martinez2018skew}, and since then many other papers analyzed their structure \cite{martinez2019reliable,martinez2019universal,ott2021bounds}  from the vector point of view. 
In this vector framework, the most prominent family of MSRD codes is given by linearized Reed-Solomon codes. They represent a sort of ``tensorization'' between Reed-Solomn codes in the Hamming metric and Gabidulin codes in the rank metric. Indeed,  these two families of codes are extremal cases of linearized Reed-Solomon codes. In our skew polynomial setting, we can define linearized Reed-Solomon codes as follows.

\begin{definition}[Linearized Reed-Solomon codes]
 The \textbf{linearized Reed-Solomon code}  of dimension $k$ is the code
 $$ \mathcal C_k^\theta\coloneqq \left\{ F \in \qspace \mid  \deg (F) <k\right\}\subseteq \qspace.$$
\end{definition}

As a trivial consequence of Theorem \ref{thm:bound_degree_improved}, we cam deduce the following well-known result.

\begin{corollary}
 Linearized Reed-Solomon codes are maximum sum-rank distance codes.
\end{corollary}

\begin{proof}
By Theorem \ref{thm:bound_degree_improved} and the definition of the metric $\dd_\Lambda$ on $\qspace$, for any nonzero skew polynomial $F\in \mathcal C_k^\theta$ we have
$$ \wt_\Lambda(F)= \ell n-\deg(\gcrd(F,H_\Lambda))\geq \ell n- \deg(F)\geq \ell n-k+1.$$
\end{proof}

\begin{remark}
 Linearized Reed-Solomon codes have been originally defined as evaluation vectors of skew polynomials, where also a derivation is considered. However, if we restrict to skew polynomials with no derivation, then we have that the codes originally defined in \cite[Definition 31]{martinez2018skew} are exactly the vectorizations of $\C_k^\theta$, that is $\C_k^\theta(\bs\alpha,\bs\mB)=\ev_{\bs\mB}(\overline{\Phi}_{\bs\alpha}(\C))$.
\end{remark}

\medskip

\subsection{Equivalence of Sum-Rank Metric Codes}
The notion of equivalence of codes in the sum-rank metric has been introduced in \cite[Theorem 2]{martinez2020hamming}. As observed there, MacWilliams extension theorem does not hold even for the case $\ell=1$ (see e.g. \cite[Example 2.9]{barra2015macwilliams}). Thus, equivalence of sum-rank metric codes is defined through isometries of the whole space. Moreover, in that paper codes were exclusively considered to be $\LL$-linear and hence only $\LL$-linear isometries were considered; see also \cite[Theorem 3.8]{alfarano2021}  In the following, we aim to characterize the  $\K$-linear and $\K$-semilinear isometries, deriving a similar result in the space $\qspace$.

The easiest setting where we can prove this result is the space $(\K^{n \times n})^\ell$.

\begin{proposition}\label{prop:klinear_equiv}
 The group of $\K$-linear isometries on $(\K^{n \times n})^\ell$ is isomorphic to $((\mathcal S_\ell \rtimes (\GL(n,\K)^\ell \times \GL(n,\K)^\ell))\rtimes (\ZZ{2})^\ell)/K$, with
 $$K=\left\{(\mathrm{id},(\lambda_1I_n,\ldots \lambda_\ell I_n),(\lambda_1^{-1}I_n,\ldots \lambda_\ell^{-1} I_n),\mathbf{0}) : \lambda_i \in \K^*\right\},$$
 and the action is given by
 $$(\pi, (M_1,\ldots,M_\ell), (N_1,\ldots,N_\ell), v) \cdot (A_1,\ldots,A_\ell) = (M_1A_{\pi^{-1}(1)}^{[v_1]}N_\ell,\ldots,M_\ell A_{\pi^{-1}(\ell)}^{[v_1]}N_\ell), $$
 where $A^{[0]}\coloneqq A$ and $A^{[1]}\coloneqq A^\top$.
\end{proposition}

\begin{proof}
 We fix the following notation. For $a\in \{1,\ldots, \ell\}$ and $M\in \K^{n\times n}$, we denote by $M^{(a)}$ the element of $(\K^{n\times n})^\ell$ that has the matrix $M$ as $a$-th entry, while all the other entries are the $0$ matrix. Moreover, for $i,j \in\{1,\ldots,n\}$ we denote by $E_{i,j}$ the matrix $e_i^\top e_j$ that has $1$ in the $(i,j)$-th entry and $0$ elsewhere. 
 
 \noindent Now, consider a $\K$-linear isometry $\phi$ of $(\K^{n\times n})^\ell$ and fix $a\in \{1,\ldots, \ell\}$. Since any matrix $E_{i,j}^{(a)}$ has sum-rank weight equal to $1$, then we have $\phi(E_{1,1}^{(a)})=M^{(b)}$, for some $b\in\{1,\ldots,\ell\}$ and $M\in\K^{n\times n}$ of rank $1$. Now, fix $j\in\{1,\ldots,n\}$ there are also  $c\in\{1,\ldots,\ell\}$ and $N\in\K^{n\times n}$ of rank $1$ such that $\phi(E_{1,j}^{(a)})=N^{(c)}$. Hence by $\K$-linearity of $\phi$ we have $$\phi((E_{1,1}+E_{1,j})^{(a)})=\phi(E_{1,1}^{(a)}+E_{1,j}^{(a)})=M^{(b)}+N^{(c)},$$
 which must have sum-rank weight equal to $1$. This is only possible if $b=c$ and $M+N$ has rank $1$. This shows that all the elementary matrices of the form $E_{1,j}^{(a)}$ are sent by $\phi$ in elements of the form $M^{(b)}$. Arguing in the same way, we can actually show that this is true for every element $E_{i,j}^{(a)}$.
 Thus, since this holds for every $a$, we can deduce that the map $\phi$ acts as
 $$ \phi(M^{(a)})=(\phi_a(M))^{(\pi^{-1}(b))},$$
 for some $\pi \in \mathcal S_\ell$, and some $\phi_{a}$ that is a $\K$-linear rank isometry of $\K^{n \times n}$. However, the $\K$-linear rank isometries are characterized and are given by the group 
 $$  (\GL(n,\K) \times \GL(n,\K))/N \rtimes \ZZ{2},$$
 where $N=\{(\lambda I_n,\lambda^{-1} I_n) : \lambda \in \K^*\}$, acting as 
 $$ (M,N,t) \cdot A = MA^{[t]}N;$$
 see e.g. \cite[Proposition 4]{morrison2014equivalence}.
 This concludes the proof.
\end{proof}

Now, using the isomorphism $\Ext_{\bs\mB}\circ\evB\circ\overline{\Phi}_{\bs\alpha}:\qspace \rightarrow (\K^{n \times n})^\ell$ and its inverse, and combining it with the characterization of the $\K$-linear isometries given in Proposition \ref{prop:klinear_equiv}, we can derive the representation of the $\K$-linear sum-rank isometries in $\qspace$.

\begin{theorem}\label{thm:equivalence_sumrank}
The group of $\K$-linear isometries on the space $\qspace$ is isomorphic to $(\mathcal S_\ell \rtimes (\GL(n,q)^\ell \times \GL(n,q)^\ell))\rtimes (\ZZ{2})^\ell$. Moreover, the action is given by
$$ F(X) \cdot (\overline{\Phi}_{\bs\alpha}^{-1}\circ \Upsilon_v \circ \overline{\Phi}_{\pi(\bs\alpha)})(\mC)\cdot G(X) $$
for some skew polynomials $F(X),G(X) \in \qspace$ such that $\wt_\Lambda(F)=\wt_\Lambda(G)=\ell n$, $\pi \in \mathcal S_\ell$ and $v\in(\ZZ{2})^\ell$.
\end{theorem}

\section{Duality of Sum-Rank Metric Codes}\label{sec:duality}

In this section we analyze a duality theory for sum-rank metric codes in the setting we proposed. We will show that the duality behaves as one would expect, and that it can be expressed in terms of the duality of the constituent skew algebras $\LL[\theta]$. This in turn can be equivalently formulated with respect to   the standard inner product duality in the vector setting and in terms of the Delsarte duality of the corresponding matrix spaces, provided that we choose  suitable $\K$-bases of $\LL$. This also implies that the dual of an MSRD code is itself MSRD. Furthermore, we analyze the dual of a linearized Reed-Solomon codes in this setting, which  turns out to be naturally (equivalent to) a linearized Reed-Solomon code.

\medskip

\subsection{Definition and Comparison with Other Dualities}
We fix the following notation and setting, that will be kept for the whole section. We choose a set of elements $\mathrm{A}=\{\alpha_1,\ldots,\alpha_\ell\}\subseteq \LL^*$ whose norms are pairwise distinct. The set $\Lambda=\Norm(\mathrm{A})=\{\lambda_1,\ldots,\lambda_\ell\}$ of norms of the elements $\alpha_i$'s will be taken to be a finite -- and hence cyclic\footnote{It is well-known that every finite multiplicative subgroup of the units of a field is cyclic} -- \textbf{subgroup} of $\K^*$. Moreover, given the vector $\bs\alpha=(\alpha_1,\ldots,\alpha_\ell)$, we will write $\bs\alpha^{-1}$ to indicate the vector whose entries are the inverses of the entries of $\bs\alpha$, that is $\bs\alpha^{-1}\coloneqq (\alpha_1^{-1},\ldots, \alpha_{\ell}^{-1})$. Observe that, if we have another set $\mathrm{A}'\coloneqq \{\beta_1,\ldots,\beta_{\ell}\}\subseteq \LL^*$  such that $\Norm(\mathrm{A}')=\Norm(\mathrm{A})=\Lambda$, then  $\overline{\Phi}_{\bs\beta}$ is a well-defined map $:\qspace \rightarrow (\LL[\theta])^\ell)$. In particular, since $\Lambda$ is a subgroup, we have that $\Norm(\{\alpha_1^{-1},\ldots,\alpha_\ell^{-1}\})=\Norm(\mathrm{A})=\Lambda$, and 
$$\overline{\Phi}_{\bs\alpha^{-1}}:\qspace \longrightarrow (\LL[\theta])^\ell $$
is well defined.

We introduce the $\K$-bilinear form $\langle \cdot,\cdot \rangle_{\Lambda}$ on $\qspace$, given by
\begin{equation}\label{eq:bilinear_form} \langle F, G\rangle_{\Lambda}=\Trace\bigg(\sum_{i=0}^{\ell n-1} f_ig_i\bigg).\end{equation}
It is easy to see that it is nondegenerate, and hence defines a duality isomorphism. Moreover, it is induced by the $\LL$-bilinear form
$$ \langle F, G\rangle_{\Lambda,\LL}=\sum_{i=0}^{\ell n-1} f_ig_i.$$

\begin{definition}
 Let $\C \subseteq \qspace$ be a (not necessarily linear) sum-rank metric code. The \textbf{dual code} of $\C$ is 
 $$\C^\perp\coloneqq \left\{ G\in \qspace \mid \langle F,G\rangle_\Lambda=0 \mbox{ for every } F \in \C \right\}.$$
 \end{definition}

Our aim is now to show that this definition of duality induced by \eqref{eq:bilinear_form} is consistent with the definitions of dualities studied in literature: the one on vector codes defined by the standard inner product, and the one of matrix codes induced by the Delsarte trace product.
As an intermediate step, we first want to contextualize the duality induced by \eqref{eq:bilinear_form} and express it in terms of the duality theory of the single components in $(\LL[\theta])^\ell$, induced by the map $\overline{\Phi}_{\bs \alpha}$.

  We first recall an auxiliary well-known result, whose proof is omitted.

\begin{lemma}\label{lem:cyclic_sum}
Let $\mathbb F$ be a field, let $\Lambda$ be a cyclic subgroup of $\mathbb F^*$ of order $\ell\geq 2$ and let $i \in \Z$. Then
$$ \sum_{\lambda \in \Lambda} \lambda ^i = \begin{cases} \ell & \mbox{ if } i \equiv 0 \mod \ell \\ 0 & \mbox{ if } i \not \equiv 0 \mod \ell.\end{cases} $$
\end{lemma}

We are now ready to express the duality of \eqref{eq:bilinear_form} in terms of the bilinear form $\langle \cdot,\cdot\rangle_{\rk}$.

\begin{proposition}\label{prop:explicit_bilinear_SR_R}
Let $\mathrm{A}=\{\alpha_1,\ldots,\alpha_\ell\}\subseteq \LL^*$ be a set of elements whose norms are pairwise distinct.  Assume that $\Lambda=\{\Norm(\alpha_1),\ldots,\Norm(\alpha_\ell)\}$ is  a finite subgroup of $\K^*$. Then
\begin{align}\label{eq:equality_bilinearform} \ell \cdot \langle F, G\rangle _{\Lambda}=\langle \overline{\Phi}_{\bs\alpha}(F),\overline{\Phi}_{\bs\alpha^{-1}}(G)\rangle_{\srk}. 
 \end{align}
\end{proposition}

\begin{proof}
First of all, we notice that, due to the $\K$-linearity of $\Trace$, it is enough to prove that 
$$ \ell \cdot \langle F,G\rangle_{\Lambda,\LL}=\sum_{\alpha \in \mathrm{A}} \langle \overline{\Phi}(F_\alpha), \overline{\Phi}(G_{\alpha^{-1}}) \rangle_{\rk,\LL}$$
Let us fix the following notation. Let $F(X)=f_0+f_1X+\ldots+f_{\ell n-1}X^{\ell n-1}$, $G(X)=g_0+g_1X+\ldots+g_{\ell n-1}X^{\ell n-1}$. Furthermore,  for any $\alpha \in \mathrm{A}$, we write 
$$\overline{\Phi}(F_\alpha)=\sum_{i=0}^{n-1}f_{\alpha,i}\theta^i, \quad \overline{\Phi}(G_\alpha)=\sum_{i=0}^{n-1}g_{\alpha,i}\theta^i. $$
Observe that, by definition we have 
\begin{equation}\label{eq:auxiliary}f_{\alpha,i}=\sum_{t=0}^{\ell-1}\NN{\theta}{i+tn}(\alpha)f_{i+t n}, \quad g_{\alpha,i}=\sum_{t=0}^{\ell-1}\NN{\theta}{i+tn}(\alpha)g_{i+t n}.\end{equation}
If we compute now the right-hand side of \eqref{eq:equality_bilinearform}, and substitute \eqref{eq:auxiliary} in it, we get
\begin{align}\label{eq:explicit_bilinear_rank}\langle \overline{\Phi}(F_\alpha), \overline{\Phi}(G_{\alpha^{-1}}) \rangle_{\rk,\LL}&=\sum_{i=0}^{n-1} f_{\alpha, i}g_{\alpha^{-1},i}=\sum_{i=0}^{n-1} \bigg(\sum_{t=0}^{\ell-1} \NN{\theta}{i+tn}(\alpha)f_{i+tn}\bigg)\bigg(\sum_{t=0}^{\ell-1} \NN{\theta}{i+tn}(\alpha^{-1})g_{i+tn}\bigg) \nonumber \\
&=\sum_{i=0}^{n-1} \bigg(\sum_{t=0}^{\ell-1}\lambda_{\alpha}^t f_{i+tn}\bigg)\bigg(\sum_{t=0}^{\ell-1} \lambda_{\alpha}^{-t}g_{i+tn}\bigg)\end{align}
Now, both the hand sides of \eqref{eq:equality_bilinearform} are $\LL$-bilinear forms, so we only need to check that they agree on the basis $\{1,X,\ldots, X^{\ell n-1}\}$. Therefore, we aim to show that
$$ \sum_{\alpha \in \mathrm{A}} \langle \overline{\Phi}((\alpha X)^{u}), \overline{\Phi}((\alpha^{-1} X)^{v}) \rangle_{\rk,\LL}=\ell \cdot \delta_{u,v}=\begin{cases} \ell & \mbox{ if } u=v \\ 0 & \mbox{ if } u \neq v. \end{cases}
$$
Let us take $F(X)=X^{i+tn}$ and $G(X)=X^{j+sn}$ for some $0\leq s,t\leq \ell -1$ and $0\leq i,j \leq n-1$. It is easy to prove that substituting these monomials $F(X)$ and $G(X)$  in \eqref{eq:explicit_bilinear_rank} we obtain $0$ whenever $i \neq j$. Therefore, let us assume that $i=j$. In this case, using \eqref{eq:explicit_bilinear_rank}, we obtain
$$ \sum_{\alpha \in \mathrm{A}} \langle \overline{\Phi}((\alpha X)^{u}), \overline{\Phi}((\alpha^{-1} X)^{v})\rangle_{\rk,\LL}=\sum_{\alpha \in \mathrm{A}} \lambda_{\alpha}^{t-s}. $$
At this point, observe that $\Lambda$ is a finite subgroup of $\K^*$, hence it is cyclic. By Lemma \ref{lem:cyclic_sum}, we have
$$\sum_{\lambda \in \Lambda}\lambda_{\alpha}^{t-s}=\begin{cases} \ell & \mbox{ if } t-s \equiv 0 \mod \ell \\ 0 & \mbox{ if } t-s \not\equiv 0 \mod \ell.
\end{cases}$$ Due to the range of $t$ and $s$, we have  $t-s \equiv 0 \mod \ell$ if and only if $t = s$. Thus, putting together all the cases, we obtain the desired result.
\end{proof}

\begin{remark}\label{rem:dualities}
 The space $\qspace$ is an $\LL$-vector space. If we write the bilinear form $\langle \cdot, \cdot\rangle_{\Lambda,\LL}$ in coordinates with respect to the monomial basis $\{1,X,\ldots, X^{\ell n-1}\}$, the matrix representing $\langle\cdot,\cdot \rangle_{\Lambda,\LL}$ is clearly the $\ell n \times \ell n$ identity matrix. On the other hand, one can consider the usual duality theory 
 inherited from  the bilinear form $\langle \cdot,\cdot \rangle_{\srk,\LL}$ over $\LL[\theta]$, i.e. 
 $$\langle F, G \rangle_{\bs\alpha,\LL}\coloneqq \langle \overline{\Phi}_{\bs\alpha}(F),\overline{\Phi}_{\bs\alpha}(G) \rangle_{\srk,\LL}. $$
 One can easily see that, if  $F(X)=X^{u}$, $G(X)=X^v$ with $u \not\equiv v \mod n$, then  $\langle F,G\rangle_{\bs\alpha,\LL}=0$.
Moreover, if $F(X)=X^{i+sn}$ and $G(X)=X^{i+tn}$ for some $0\leq s,t\leq \ell -1$ and $0\leq i \leq n-1$,  using \eqref{eq:auxiliary} it is easy  to derive that 
$$\langle X^{i+sn},X^{i+tn}\rangle_{\bs\alpha,\LL}=\sum_{\alpha \in \mathrm{A}} \NN{\theta}{i}(\alpha^2)\lambda_\alpha^{s+t}.$$
Thus, the matrix associated to $\langle \cdot,\cdot\rangle_{\bs\alpha,\LL}$ with respect to the monomial basis is given by
$$\sum_{\alpha \in \mathrm{A}} M_{\alpha}, $$
where
\begin{align*}M_{\alpha}&\coloneqq (\lambda_{\alpha}^{s+t})_{s,t}\otimes \diag(\NN{\theta}{0}(\alpha^2),\ldots,\NN{\theta}{n-1}(\alpha^2)) \in \LL^{\ell n \times \ell n}.
\end{align*}

One could directly consider the canonical orthonormal basis of $(\LL[\theta])^\ell$ for the $\LL$-bilinear form $\langle\cdot , \cdot \rangle_{\srk,\LL}$, namely 
$$ \{(\theta^{i}, 0,\ldots,0) \mid 0\leq i < n\}\cup\{(0,\theta^{i},\ldots,0) \mid 0\leq i < n\} \cup \ldots \cup \{(0,0,\ldots,\theta^{i}) \mid 0\leq i < n\}.  $$
It is easy to see that this is the image under the map $\overline{\Phi}_{\bs\alpha}$ of the set of skew polynomials
$$ \mathcal M\coloneqq \left\{\ell^{-1}\NN{\theta}{i}(\alpha_j^{-1}){\lambda_j}(X^{\ell n}-1)(X^n-\lambda_j)^{-1}X^i \mid 1\leq j \leq \ell, \, 0 \leq i \leq n-1\right\}. $$
Indeed, 
\begin{align*}((X^{\ell n}-1)(X^n-\lambda_j)^{-1}X^i)_{\alpha_t}&=((\alpha_tX)^{\ell n}-1)((\alpha_t X)^n-\lambda_j)^{-1}(\alpha_t X)^{i}\\ 
&=\lambda_t^{-1}(X^{\ell n}-1)(X-\lambda_t^{-1}\lambda_j)^{-1}\NN{\theta}{i}(\alpha_t)X^i,
\end{align*}
which is $0$ modulo $(X^n-1)$ whenever $t\neq j$. If $t=j$, we have
\begin{align*}((X^{\ell n}-1)(X^n-\lambda_j)^{-1}X^i)_{\alpha_j}
&=\lambda_j^{-1}(X^{\ell n}-1)(X-1)^{-1}\NN{\theta}{i}(\alpha_j)X^i \\
&=\lambda_j^{-1}\NN{\theta}{i}(\alpha_j)(1+X^n+\ldots+X^{(\ell-1)n}) X^i,
\end{align*}
that equals to $\lambda_j^{-1}\NN{\theta}{i}(\alpha_j)\ell X^i$ modulo $(X^n-1)$. Thus,  one also obtains that
$$\sum_{\alpha\in \mathrm{A}}M_\alpha = N^\top N,$$
where $N$ is the change-of-basis matrix from $\mathcal M$ to $\mB=\{1,X,\ldots,X^{\ell n-1}\}$. Observe that the division by $\ell$ is always well-defined, since $\ell$ is the order of a cyclic subgroup of $\K^*$, and hence it cannot be a multiple of  $\ch(\K)$.
\end{remark}

Now, we want to analyze the connection with the duality of codes in $\LL^{\ell n}$ with respect to the standard inner product. In this setting, if we have a code $C\subseteq \LL^{\ell n}$, we define the dual code, to be the orthogonal space $C^{\perp_{\vv}}$ with respect to the bilinear form
$$ \langle \mathbf{u},\mathbf{v} \rangle_{\vv} := \Trace(\mathbf{u}\mathbf{v}^\top), \qquad \mbox{ for every } \mathbf{u},\mathbf{v} \in \LL^{\ell n}.$$
First, for a given $\K$-basis $\mB=(\beta_1,\ldots,\beta_n)$ of $\LL$, we denote by $\mB^*:=(\beta_1^*,\ldots,\beta_n^*)$ its dual basis with respect to the trace bilinear form induced by $\Trace$, that is such that
$$ \Trace(\beta_i\beta_j^*)=\delta_{i,j}=\begin{cases} 1 & \mbox{ if } i=j, \\
0 & \mbox{ if } i \neq j.\end{cases}$$
\begin{theorem}\label{thm:duality_vector_poly}
 Assume that $\Lambda=\Norm(\{\alpha_1,\ldots,\alpha_\ell\})$ is a multiplicative subgroup of $\K^*$, and let $\bs\mB:=(\mB_1,\ldots,\mB_{\ell})$ be a vector of $\K$-bases of $\LL$. Let $\C\subseteq \qspace$ be a $\K$-linear code. Then
 $$ \C^\perp(\bs\alpha, \bs\mB)=(\C(\bs\alpha^{-1},\bs\mB^*))^{\perp_{\vv}},$$
 where $\bs\mB^*:=(\mB_1^*,\ldots,\mB_\ell^*)$.
\end{theorem}

\begin{proof}
 Observe that it is enough to verify that 
 $$ ((\evB \circ \overline{\Phi}_{\bs\alpha})(X^u))( (\ev_{\bs\mB^*} \circ \overline{\Phi}_{\bs\alpha^{-1}})(X^v))^\top= \omega \delta_{u,v},$$
 for some $\omega \in \K^*$. 
  Let us fix $u,v \in \{0,\ldots, \ell n-1\}$ and let us write $u=i+tn$ and $v=j+sn$, for some $0 \leq s,t \leq \ell -1$ and $0 \leq i,j,\leq n-1$.
  We can write
 \begin{equation}\label{eq:duality_vector} ((\evB \circ \overline{\Phi}_{\bs\alpha})(X^u))( (\ev_{\bs\mB^*} \circ \overline{\Phi}_{\bs\alpha^{-1}})(X^v))^\top=\sum_{i=1}^\ell (\ev_{\mB_i}(\overline{\Phi}((\alpha_iX)^u))(\ev_{\mB_i^*}(\overline{\Phi}((\alpha_i^{-1}X)^v)))^\top.
 \end{equation}
 We now distinguish two cases.
 
 \noindent \underline{\textbf{Case I: Normal bases.}}
 Let $\alpha \in \LL^*$ such that $\Norm(\alpha)=\lambda\in\Lambda$ and let us fix a normal basis $\mE=(\beta,\theta(\beta),\ldots,\theta^{n-1}(\beta))$. In this case, it is known that $\mE^*$ is also a normal basis, i.e. $\mE^*=(\gamma,\theta(\gamma),\ldots,\theta^{n-1}(\gamma))$. Then we have
 \begin{align}\label{eq:duality_vector_component}(\ev_{\mE}(\overline{\Phi}((\alpha X)^u))(\ev_{\mE^*}(\overline{\Phi}((\alpha^{-1}X)^v)))^\top&=(\ev_{\mE}(\lambda^t\NN{\theta}{i}(\alpha)\theta^i))(\ev_{\mE^*}(\lambda^{-s}\NN{\theta}{j}(\alpha^{-1})\theta^j))^\top \nonumber\\
 &=\lambda^{t-s}\NN{\theta}{i}(\alpha)\NN{\theta}{j}(\alpha^{-1})\ev_{\mE}(\theta^i)\ev_{\mE^*}(\theta^j) \nonumber\\ 
 &=\lambda^{t-s}\NN{\theta}{i}(\alpha)\NN{\theta}{j}(\alpha^{-1})\sum_{a=0}^{n-1}\theta^{i+a}(\beta)\theta^{j+a}(\gamma) \nonumber\\
 &= \lambda^{t-s}\NN{\theta}{i}(\alpha)\NN{\theta}{j}(\alpha^{-1})\Trace(\theta^{i}(\beta)\theta^j(\gamma)) \nonumber\\
 &=\lambda^{t-s}\NN{\theta}{i}(\alpha)\NN{\theta}{j}(\alpha^{-1})\delta_{i,j} \nonumber \\
 &=\lambda^{t-s}\delta_{i,j},
 \end{align}
 where the last equality follows from the fact computed quantity is nonzero only  when $i=j$, in which case $\NN{\theta}{i}(\alpha)\NN{\theta}{j}(\alpha^{-1})=1$. Assume now  that each basis $\mB_i$ is normal. Then, combining \eqref{eq:duality_vector} and \eqref{eq:duality_vector_component}, we obtain
 $$ ((\evB \circ \overline{\Phi}_{\bs\alpha})(X^u))( (\ev_{\bs\mB^*} \circ \overline{\Phi}_{\bs\alpha^{-1}})(X^v))^\top=\sum_{i=1}^\ell \lambda^{t-s}\delta_{i,j}=\ell \delta_{i,j}\delta_{t,s}=\ell \delta_{u,v},$$
 where the second to last equality follows from Lemma \ref{lem:cyclic_sum}.
 
  \noindent \underline{\textbf{Case II: General bases.}}
    Now, suppose that $\mB$ is a general
  $\K$-basis of $\LL$. There exists $X\in\GL(n,\K)$ such that $\mB=\mE
  X$, where $\mE$ is normal. Moreover, we also have that $\mB^*=\mE^* (X^{-1})^\top$. Hence, we
  get
  \begin{align*}
    \ev_\mB(\theta^i) \ev_{\mB^*}(\theta^j)^\top &=  \ev_{\mE X}(\theta^) \ev_{\mE^* (X^{-1})^\top}(\theta^j)^\top \\
    &= \ev_\mE(\theta^i)X(\ev_{\mE^*}(\theta^j)(X^{-1})^\top )^\top \\
    &=\ev_\mE(\theta^i)\ev_{\mE^*}(\theta^j)^\top, 
  \end{align*}
  and we can conclude using Case I.
\end{proof}

Now, we proceed relating the duality between vector sum-rank metric codes and matrix sum-rank metric codes. By transitivity, we then also deduce a relation of the duality between the skew polynomial and the matrix frameworks. 

In the following if we have a code $C\subseteq (\K^{n\times n})^\ell$, we define the dual code, to be the orthogonal space $C^{\perp_{\MM}}$ with respect to the bilinear form
$$ \langle (M_1,\ldots,M_\ell),(N_1,\ldots,N_\ell) \rangle_{\MM} := \sum_{i=1}^\ell\Tr(M_iN_i^\top), \qquad \mbox{ for every } M_i,N_i \in \K^{n\times n}.$$
We omit the proof of the following result, since it can be derived in exactly the same way as done for the rank metric; see e.g. \cite{grant2008duality,ravagnani2016rank}.
\begin{proposition}\label{prop:duality_vector_matrix}
 Let $C\subseteq \LL^{\ell n}$ be a $\K$-linear sum-rank metric code. Then, for every $\ell$-uple of $\K$-bases of $\LL$ we have
 $$ \Ext_{\bs\mB}(C^{\perp_{\vv}})=\Ext_{\bs\mB^*}(C)^{\perp_{\MM}}.$$
\end{proposition}

\medskip

\subsection{Duality of MSRD and Linearized Reed-Solomon Codes}

As a consequence, we can immediately deduce a duality result for MSRD codes, that holds also in this skew polynomials framework.

\begin{theorem}\label{thm:dualMSRD}
 Let $\C\subseteq \qspace$ be a $\K$-linear MSRD code and assume that $\Lambda$ is a finite group. Then $\C^\perp$ is a $\K$-linear MSRD code.  
\end{theorem}

\begin{proof}
 Let $\C\subseteq \qspace$ be a $\K$-linear MSRD code. For any vector $\bs\alpha=(\alpha_1,\ldots,\alpha_\ell)\in\LL^\ell$ such that $\Norm(\{\alpha_1,\ldots,\alpha_\ell\})=\Lambda$, and for every $\ell$-uples  $\bs\mB=(\mB_1,\ldots,\mB_\ell)$ $\bs\mE=(\mE_1,\ldots,\mE_\ell)$  of $\K$-bases of $\LL$, the code $(\Ext_{\bs\mE}\circ \evB\circ \overline{\Phi}_{\bs\alpha})(\C)=\Mat_{\bs\alpha,\bs\mB,\bs\mE}(\C)\subseteq (\K^{n \times n})^\ell$ is also MSRD. This is due to the fact that the map 
 $$\Ext_{\bs\mE}\circ \evB\circ \overline{\Phi}_{\bs\alpha}:(\qspace,\dd_{\Lambda})\longrightarrow((\K^{n \times n})^\ell,\dd_{\MM})$$
 is a $\K$-linear isometry. Furthermore, the code $\Mat_{\bs\alpha,\bs\mB,\bs\mE}(\C)^{\perp_{\MM}}$ is also MSRD, due to \cite[Theorem VI.1]{byrne2020fundamental}.\footnote{The  result in \cite{byrne2020fundamental} is only stated for $\K=\Fq$, but  the proof can be adapted straightforwardly to any field $\K$.} Now using Proposition \ref{prop:duality_vector_matrix} and Theorem \ref{thm:duality_vector_poly}, we get
 \begin{align*} \Mat_{\bs\alpha,\bs\mB,\bs\mE}(\C)^{\perp_{\MM}}&=(\Ext_{\bs\mE}(\C(\bs\alpha,\bs\mB)))^{\perp_{\MM}}=\Ext_{\mE^*}(\C(\bs\alpha,\bs\mB)^{\perp_{\vv}})\\ &=\Ext_{\bs\mE^*}(\C^\perp(\bs\alpha^{-1},\mB^*))=\Mat_{\bs\alpha^{-1},\bs\mB^*,\bs\mE^*}(\C^\perp).
 \end{align*}
 Since also the map 
  $$\Ext_{\bs\mE^*}\circ \ev_{\mB^*}\circ \overline{\Phi}_{\bs\alpha^{-1}}:(\qspace,\dd_{\Lambda})\longrightarrow((\K^{n \times n})^\ell,\dd_{\MM})$$
  is a $\K$-linear isometry, we can conclude that $\C^\perp$ is MSRD.
\end{proof}

\begin{remark}
 In the proof of Theorem \ref{thm:dualMSRD} we relied on the duality statement of \cite[Theorem VI.1]{byrne2020fundamental} in the matrix setting, although the duality of MSRD codes was already proved in \cite[Theorem 5]{martinez2019theory} in the vector framework. However,  the latter result was only given for $\LL$-linear codes, while the one in \cite{byrne2020fundamental} was extended to $\K$-linear sum-rank metric codes.
\end{remark}

As observed in Remark \ref{rem:dualities}, it is clear that a duality theory based on the bilinear form $\langle \cdot,\cdot\rangle_{\srk}$ is not ``monomial-friendly'' -- that is, the monomial basis is not orthonormal with respect to $\langle \cdot,\cdot\rangle_{\srk,\LL}$ -- and this makes the computations more complicated. The bilinear form $\langle \cdot,\cdot\rangle_{\Lambda}$ is indeed more natural in this setting, as it allows to show straightforwardly that the dual of a Linearized Reed-Solomon code is equivalent to a  Linearized Reed-Solomon code.

\begin{theorem}
 Let $k\leq \ell n$ be a positive integer. Then 
 $$(\C_k^\theta)^\perp=\C_{\ell n-k}^\theta \cdot X^{k}=X^k\cdot \C_{\ell n-k}^\theta.$$   
\end{theorem}

\begin{proof}
This is  immediate, since $(\C_k^\theta)^\perp=\langle \{X^{i} \colon k\leq i \leq \ell n -1\}\rangle_\LL$.
\end{proof}

\begin{remark}
In \cite[Theorem 4]{martinez2019reliable} it was shown that the dual of a Linearized Reed-Solomon code over a finite field is in turn a Linearized Reed-Solomon code, in the vector framework. The hypothesis of working over a finite field could not be removed, since the proof relies on the choice of a suitable primitive element of $\LL=\F_{q^n}$. Our approach is more intrinsic and natural, and allows to show the same result \emph{over any cyclic extension of fields}. Moreover, since $\wt_{\Lambda}(X^k)=\ell n$, then, by Theorem \ref{thm:equivalence_sumrank}, $(\C_k^\theta)^\perp$ is equivalent to  $\C_{\ell n-k}^\theta$. In order to get the result on the vector representation as the one in \cite[Theorem 4]{martinez2019reliable}, at this point one can simply use Theorem \ref{thm:duality_vector_poly}, obtaining that $\C_k^\theta(\bs\alpha,\bs\mB)^{\perp_{\vv}}=\mathcal D(\bs\alpha^{-1},\bs\mB^*)$, where $\mathcal D=\C_{\ell n-k}^\theta\cdot X^{k}$. It is not difficult to check that $\mathcal D(\bs\alpha^{-1},\bs\mB^*)=\C_{\ell n-k}^\theta(\bs\alpha^{-1},\bs\mE)$, where $\bs\mE=(\mE_1,\ldots,\mE_\ell)$ and $\mE_i=\NN{\theta}{k}(\alpha_i^{-1})\cdot\theta^k(\mB_i^*)$ for every $i\in\{1,\ldots,\ell\}$.
\end{remark}

The following is a very easy calculation showing that the adjoint of a linearized Reed-Solomon code is equivalent to a linearized Reed-Solomon code.
\begin{proposition}
 Assume that $\Lambda$ is a subgroup of $\K^*$, and let $k\leq \ell n$ be a positive integer. Then
 $$ (\C_k^\theta)^\top= \C_k^\theta\cdot X^{\ell n -k+1}=X^{\ell n -k+1}\cdot\C_k^\theta.$$
\end{proposition}

\section{Twisted Linearized Reed-Solomon Codes}\label{sec:twisted_LRS}

In this section we introduce a twisted version of linearized Reed-Solomon codes in the sum-rank metric, which generalizes the notion of twisted Gabidulin codes defined in \cite{sheekey2016new}. We will show that, under certain assumptions, these codes are maximum sum-rank distance codes. Moreover, we will see that when we specialize our construction to only one block, this results in a twisted Gabidulin code. Note that for the moment we are still considering a general cyclic Galois extension $\LL/\K$.

\medskip

\subsection{Definition and Properties}

We start with a preliminary result that generalizes \cite[Lemma 3]{sheekey2016new}.
\begin{proposition}\label{prop:full_kernel_norm_condition}
   Let $\Lambda\subseteq \K^*$ and let $F(X)=f_0 +\ldots +f_kX^{k}$  be such that  $k=\deg (F)=\sum_{\lambda\in \Lambda}d_\lambda(F)$. Then 
   $$\Nn_{\LL/\K}(f_0/f_k)=(-1)^{kn}\prod_{\lambda \in \Lambda}\lambda^{d_\lambda(F)}.$$ 
\end{proposition}

\begin{proof}
 Assume that  $k=\sum_{\lambda \in \Lambda}d_\lambda(F)=\sum_{\lambda\in \Lambda}\dim_{\LL}(\ker(A_F-\lambda I_k))$.  Therefore, by Theorem \ref{thm:bound_degree_improved}, the set $\overline{\Lambda}\coloneqq \{\lambda \in \Lambda \mid d_{\lambda}(F)>0\}$ is the set of eigenvalues for $A_F$, where each eigenvalue $\lambda$ has multiplicity $d_\lambda(F)$. This implies that $$\det(A_F)=\prod_{\lambda\in \overline{\Lambda}}\lambda^{d_\lambda(F)}.$$
 On the other hand, by definition of the matrix $A_F$,
 $$\det(A_F)=\prod_{i=0}^{n-1}\det(\theta^i(C_F))=\prod_{i=0}^{n-1}\theta^i(\det(C_F))=\Nn_{\LL/\K}((-1)^kf_0/f_k)=(-1)^{kn}\Nn_{\LL/\K}(f_0/f_k).$$
 Combining the two equalities we obtain the desired result.
\end{proof}

\begin{definition}\label{def:twistedLRS}
 Let $k,n, \ell$ be positive integers with $1\leq k \leq \ell n$. Let  $\Lambda=\{\lambda_1,\ldots,\lambda_\ell\}\subseteq \K^*$ be a finite set. 
 Furthermore, let $\eta \in \LL$ such that $(-1)^{kn}\Nn_{\LL/\K}(\eta) \not\in \langle \Lambda\rangle$, where $\langle \Lambda\rangle$ denotes the multiplicative subgroup of $\K^*$ generated by $\Lambda$.  The code 
  $$ \mathcal L_k^{\theta}(\eta,h)\coloneqq \left\{f_0+ \ldots+ f_{k-1}X^{k-1} +\eta \theta^h(f_0)X^{k} \mid f_i \in \LL  \right\}\subseteq \qspace $$
 is called a \textbf{twisted linearized Reed-Solomon code}.
\end{definition}

\begin{theorem}\label{thm:twisted_are_MSR}
 The code $\mathcal L_k^{\theta}(\eta,h)$ is a maximum sum-rank distance code.
\end{theorem}

\begin{proof}
First, notice that $\mathcal L_k^{\theta}(\eta,h)$ is $\K$-linear, with $\dim_{\K}(\mathcal L_k^{\theta}(\eta,h))=nk=[\LL:\K]k$. Therefore, by Theorem \ref{thm:singleton} we need to prove that for every $F\in \mathcal L_k^{\theta}(\eta,h)$ we have $\wt_\Lambda(F)\geq \ell n-k+1$. Observe that, by definition of $F_{\bs\alpha}$,  we have
$$ \wt_{\Lambda}(F)= \sum_{i=1}^\ell(n-\dim_{\K}(\ker(F_{\alpha_i}))= \ell n-\sum_{i=1}^\ell d_{\lambda_i}(F),$$
 where the last equality comes from Proposition \ref{prop:equality_d_lambda}.
 If $\deg (F)\leq k-1$, then by Theorem \ref{thm:bound_degree_improved}, we have the desired inequality. Hence, assume that $\deg (F)=k$. If $\sum_{i=1}^\ell d_{\lambda_i}(F)=k$, then by Proposition \ref{prop:full_kernel_norm_condition} we must have $$(-1)^{kn}\Nn_{\LL/\K}(\eta)=(-1)^{kn}\Nn_{\LL/\K}(f_0/f_k)=\prod_{i=1}^\ell \lambda_i^{d_{\lambda_i}(F)},$$ which contradicts the assumption on $\eta$.
\end{proof}

\begin{remark}
If we set $\ell=1$ and choose $\bs\alpha=\alpha_1=1$, then  $\mL_k^\theta(\eta,h)$ is simply the twisted Gabidulin code $\mH_k^\theta(\eta,h)$.  Moreover, we can see that if we choose $\eta= 0$, then $\mL_k^\theta(0,h)=\C_k^\theta$ is a linearized Reed-Solomon code. For this reason, we will refer to \emph{nontrivial twisted linearized Reed-Solomon codes} to denote those with $\eta \neq 0$.
\end{remark}

\begin{remark}
 If we set $\LL=\K$ -- or equivalently $\theta=\mathrm{id}$ -- the sum-rank metric is in fact just the Hamming metric. Moreover, in this case, the code $\mL_k^{\mathrm{id}}(\eta,h)=\mL_k^{\mathrm{id}}(\eta,0)$ coincides with the twisted Reed-Solomon code proposed in \cite[Definition 7]{beelen2017twisted}. Also the condition on $\eta$ in Definition \ref{def:twistedLRS} turns out to be the same as the one proposed there. See also \cite{beelen2018structural}.
\end{remark}

 We now illustrate how to represent twisted linearized Reed-Solomon codes also in the vector and matrix frameworks with the aid of the following example.
 
\begin{example}
 Let us fix the same setting used in Examples \ref{exa:f5}, \ref{exa:f5new} and \ref{exa:f5bis}, that is $\K=\F_5$, $\LL=\F_{5^3}$ with $\gamma$ primitive element of $\LL$ which is a root of $y^3+3y+3$. However, this time we choose $\ell=2$ with $\bs\alpha=(1,4)$, and hence $\Lambda=\{1,4\}$. Let us take a twisted linearized Reed-Solomon code which is linear over $\LL=\F_{5^3}$, that is with $h=0$. Observe that $\Norm(2)=3$ and consider the code
 $$ \mC:=\mL_2^\theta(2,0)=\{f_0+f_1X+2f_2X^2 \mid f_0,f_1 \in \F_{5^3} \}=\langle X,1+2X^2\rangle_{\F_{5^3}}.$$
 By Theorem \ref{thm:twisted_are_MSR}, $\mL_2^\theta(2,0)$ is a MSRD code, with minimum sum-rank distance
 $\dd_\Lambda(\mC)=5$. Let $\mB=(1,\gamma, \gamma^2)$ and define $\bs\mB=(\mB,\mB)$. We now compute the vector representation with respect to $\bs\mB$ and $\bs\alpha$. First,  let us call $F(X)\coloneqq X$ and
 $G(X)\coloneqq 1+2X^2$ and compute
  $$\begin{array}{rclrcl}
    \overline{\Phi}(F_1)& \!\!\!\!=\!\!\!\!& \theta, \;\;\; &\;\;\;\overline{\Phi}(G_1)&\!\!\!\!=\!\!\!\!&2\theta^2 + \mathrm{id}, \\
    \overline{\Phi}(F_4)&\!\!\!\!=\!\!\!\!& 2\theta, & \overline{\Phi}(G_4)&\!\!\!\!=\!\!\!\!& 3\theta^2 +   \mathrm{id}.
 \end{array}$$
 Then, 
 \begin{align*} \C(\bs\alpha,\bs\mB)&=\langle \evB(\overline{\Phi}_{\bs\alpha}(F)),\evB(\overline{\Phi}_{\bs\alpha}(G)) \rangle_{\F_{5^3}}  \\
 &= \mathrm{Rowsp}\left(\begin{array}{ccc|ccc} 1 & \gamma^{5} & \gamma^{10} & \gamma^{31} & \gamma^{36} & \gamma^{41} \\
 \gamma^{93} & \gamma^{115} & \gamma^{29} & \gamma^{62} & \gamma^{54} & \gamma^{51}  \\ 
 \end{array}\right).
 \end{align*}
 Let $u$ and $v$ denote the first and the second row of the above matrix, and let us fix $\bs\mE\coloneqq \bs\mB$. 
 Then we have
 $$ \Ext_{\bs\mB}(u)=\left(\begin{pmatrix} 1 & 4 & 3 \\ 0 & 4 & 2 \\ 0 & 2 & 0\end{pmatrix},\begin{pmatrix} 2 & 3 & 1 \\ 0 & 3 & 4 \\ 0 & 4 & 0\end{pmatrix}\right), $$
  $$ \Ext_{\bs\mB}(v)=\left(\begin{pmatrix} 3 & 2 & 2 \\ 0 & 1 & 1 \\ 0 & 1 & 4\end{pmatrix},\begin{pmatrix}4 & 3 & 3 \\ 0 & 1 & 4 \\ 0 & 4 & 3 \end{pmatrix}\right). $$
  If $M$ denotes the companion matrix of the polynomial $y^3+3y+3$, that is 
  $$ M=\begin{pmatrix} 0 & 0 & 2 \\ 1 & 0 & 2 \\ 0 & 1 & 0
  \end{pmatrix},$$
  then we have 
  $$ \Mat_{\bs\alpha,\bs\mB,\bs\mB}(\C)=\langle \{M^i \cdot  \Ext_{\bs\mB}(u) : 0 \leq i \leq 2\} \cup \{ M^i\cdot\Ext_{\bs\mB}(v) : 0 \leq i \leq 2\}\rangle_{\F_5}.$$
\end{example}

\begin{remark}
 Twisted Gabidulin codes have been generalized to a wider class of MRD codes, called \emph{additive} (\emph{generalized}) \emph{twisted Gabidulin codes}, by Otal and Ozbudak in \cite{otal2016additive}. Their main idea is that to extend the definition to any automorphism of $\LL$ acting on the coefficient of $X^k$. Formally, they considered the set of $\theta$-polynomials
 $$\left\{f_0\mathrm{id}+\ldots+f_{k-1}\theta^{k-1}+\eta \tau^h(f_0)\theta^{k} \mid f_i\in \LL \right\}\subseteq \LL[\theta],$$
 where $\tau$ is any element in $\Aut(\LL)$. If we call $u\coloneqq [\LL:\LL^{\tau}]$ and we assume that $(-1)^{uk}\Norm(\eta)\neq 1$, they showed that the resulting code is MRD. Also in our setting we can easily generalize our family of codes to a class of \textbf{additive twisted linearized Reed-Solomon codes}. If we fix $\tau \in \Aut(\LL)$ of finite order $u$, we can define the code
 $$\mL_k^\tau(\eta,h)\coloneqq \left\{f_0+ \ldots+ f_{k-1}X^{k-1} +\eta \tau^h(f_0)X^{k} \mid f_i \in \LL  \right\}\subseteq \qspace. $$
 It is easy to verify that under the assumption that $(-1)^{uk}\Norm(\eta) \notin \langle \Lambda\rangle$, the resulting  code is maximum sum-rank distance. Furthermore, such a code is $\LL^{\tau}$-linear.  
\end{remark}

\medskip

\subsection{Existence and Parameters}

In principle, one could define the  code $\mL_k^\theta(\eta,h)$ without any assumption on the value $\eta$. However, in this case it is not guaranteed that we have a maximum sum-rank distance code. Indeed, the condition 
\begin{equation}\label{eq:norm_condition}
    (-1)^{kn}\Nn_{\LL/\K}(\eta) \not\in \langle \Nn_{\LL/\K}(\alpha_1),\ldots, \Nn_{\LL/\K}(\alpha_\ell)\rangle
\end{equation} 
is used in the proof of Theorem \ref{thm:twisted_are_MSR} to show that  the code is maximum sum-rank distance. This condition could produce a restriction on $\ell$, that is the number of blocks allowed in the construction of nontrivial twisted linearized Reed-Solomon codes. This clearly depends on how big the subgroup $\langle \Lambda \rangle$ is. For instance, if we choose $\Lambda$ such that $\langle \Lambda\rangle =\K^*$, then the only element $\eta$ satisfying \eqref{eq:norm_condition} is $\eta=0$, and hence we cannot construct any nontrivial twisted linearized Reed-Solomon code of the form $\mL_k^\theta(\eta,h)$. 

Here we investigate the possible values $\ell$ and the existence of twisted linearized Reed-Solomon codes.

\begin{proposition}\label{prop:number_of_blocks_twisted}\,
\begin{enumerate}
    \item If $\K=\Fq$, then we can construct a nontrivial twisted linearized Reed-Solomon code $\mL_k^\theta(\eta,h)$ with $\ell$ blocks satisfying \eqref{eq:norm_condition}, whenever  $\ell \leq \frac{q-1}{r}$, where $r$ is the smallest prime dividing $(q-1)$.
    \item If $\K$ is an infinite field admitting a cyclic Galois extension $\LL$ of degree $n$, then we can construct a nontrivial twisted linearized Reed-Solomon code $\mL_k^\theta(\eta,h)$ with $\ell$ blocks satisfying \eqref{eq:norm_condition}, for every $\ell \in \mathbb N$.
    \end{enumerate}
\end{proposition}

\begin{proof}
According to  \eqref{eq:norm_condition}, we can   construct nontrivial twisted linearized Reed-Solomon if and only if there exists $\lambda_1,\ldots, \lambda_\ell$ distinct elements of $\K^*$ all belonging to  a proper subgroup $H$ of $\K^*$. 
\begin{enumerate}
\item In the case of finite fields, the biggest proper subgroup of $\Fq^*$ has cardinality $\frac{q-1}{r}$, where $r$ is the smallest prime dividing $q-1$.
\item If $\K$ is an infinite field, then it is well-known that $\K^*$ is not finitely generated. Hence, $\langle \Nn_{\LL/\K}(\alpha_1),\ldots, \Nn_{\LL/\K}(\alpha_\ell)\rangle\neq \K^*$, for any choice of $\ell\in \N$ and we can always find an element $\eta \in \K^*\setminus\langle \Nn_{\LL/\K}(\alpha_1),\ldots, \Nn_{\LL/\K}(\alpha_\ell)\rangle $. \end{enumerate}\end{proof}

Hence, if we compare twisted linearized Reed-Solomon codes with linearized Reed-Solomon codes, Proposition \ref{prop:number_of_blocks_twisted} shows that the price to pay in terms of number of blocks is only in the case of finite fields.

\begin{remark}
 Part (2) of Proposition \ref{prop:number_of_blocks_twisted} is consistent with the fact that there exist nontrivial twisted Gabidulin codes over a finite field $\Fq$ whenever $q\neq 2$. Indeed, as a particular case of Proposition \ref{prop:number_of_blocks_twisted}, it is  possible to find a proper subgroup of $\Fq^*$ of order $\ell=1$ if and only if $q>2$, which is the same condition on the element $\eta$ imposed in the definition of nontrivial twisted Gabidulin codes.
\end{remark}

\begin{remark}
 In the case of finite fields, one could try to construct  codes with the same shape of twisted linearized Reed-Solomon codes removing the condition \eqref{eq:norm_condition}, hoping to obtain longer MSRD codes. Indeed, this assumption is sufficient but not necessary in general for producing MSRD codes, and it can be replaced with the condition
 $$ (-1)^{kn}\Norm(\eta) \notin \bigg\{\prod_{i=1}^\ell \lambda_i^{j_i} : j_1,\ldots,j_\ell \in \N, j_1+\ldots+j_\ell=k \bigg\}.$$
 However, when $k\geq 3$ and $q\geq 13$ is odd, there is no gain in the length; see \cite[Theorem 3.1]{roth1992t}.
\end{remark}

\medskip

\subsection{Dual and Adjoint}

In this section we conclude the study on twisted linearized Reed-Solomon codes by deriving their dual and their adjoint codes.

\begin{theorem}\label{thm:dual_twisted} The dual code of a twisted linearized Reed-Solomon code is given by
 $$\mathcal L_k^\theta(\eta,h)^\perp = \mathcal L_{\ell n-k}^\theta(-\theta^{n-h}(\eta), n-h)\cdot X^k.$$
\end{theorem}

\begin{proof}
 It is immediate to see that every monomial $\lambda X^i$, with $k+1\leq i \leq \ell n-1$ belongs to  $\mathcal L_k^\theta(\eta,h)^\perp$. It is enough to show that also every binomial of the form $$F^\mu(X)\coloneqq\mu X^k-\theta^{\ell n-h}(\eta)\theta^{\ell n -h}(\mu)$$ belongs to $\mathcal L_k^\theta(\eta,h)^\perp$. Clearly $F^\mu(X)$ is orthogonal to all the monomials $\beta X^j$ for $1\leq j \leq k-1$. It remains to show that for every $\beta \in \LL$ we have
 $$ \langle F^\mu(X), G^\beta(X) \rangle_{\Lambda}=0,$$
 where $G^\beta(X)=\eta \theta^h(\beta)X^k+\beta$. An easy computation shows that
 $$ \langle F^\mu(X),G^\beta(X)\rangle_{\Lambda}=\Trace(\langle F^\mu(X),G^\beta(X)\rangle_{\Lambda,\LL})=\Trace(\mu\eta\theta^h(\beta)-\theta^{n-h}(\eta)\theta^{n-h}(\mu)\beta)=0. $$
 Thus, we have
 \begin{align*}\mathcal L_k^\theta(\eta,h)^\perp&= \left\{f_0X^k+f_{1}X^{k+1}+\ldots+f_{n-k+1}X^{\ell n-1}-\theta^{n-h}(f_{0})\theta^{n-h}(\eta) : f_i \in \LL \right\}\\ &= \mathcal L_{\ell n-k}^\theta(-\theta^{n-h}(\eta), n-h)\cdot X^k.\end{align*}
\end{proof}

\begin{theorem}\label{thm:adjoint_twisted} Let $\Lambda$ be a cyclic group. Then, 
the adjoint code of a twisted linearized Reed-Solomon code is given by
$$\mathcal L_k^\theta(\eta,h)^\top=\mL_{k}^{\theta}(\theta^{n-k}(\eta^{-1}),k-h)\cdot X^{\ell n-k}.$$
\end{theorem}

\begin{proof}
 By Theorem \ref{thm:adjoint}, the adjoint of any monomial $X^{i}$ is $(X^i)^\top=X^{\ell n-i}$. Morover, for every polynomial $G^\beta(X):=\eta \theta^h(\beta)X^k+\beta$, we can compute the adjoint using again Theorem \ref{thm:adjoint} and obtaining
 $$ (G^\beta)^\top(X)=\beta+\theta^{n-k}(\eta)\theta^{n-k+h}(\beta)X^{\ell n-k}. $$
  Thus, we have
 \begin{align*}\mathcal L_k^\theta(\eta,h)^\top&= \left\{\theta^{n-k}(\eta)\theta^{n-k+h}(f_0)X^{\ell n-k}+f_{1}X^{\ell n-k+1}+\ldots+f_{k-1}X^{\ell n-1}+f_0 : f_i \in \LL \right\}\\
 &= \left\{f_0X^{\ell n-k}+f_{1}X^{\ell n-k+1}+\ldots+f_{k-1}X^{\ell n-1}+\theta^{n-k}(\eta^{-1})\theta^{k-h}(f_0) : f_i \in \LL \right\}\\&= \mL_{k}^{\theta}(\theta^{n-k}(\eta^{-1}),k-h)\cdot X^{\ell n-k}.\end{align*}
\end{proof}

\section{The Analogue of Trombetti-Zhou Construction and New MDS Codes}\label{sec:TZ}

Another family of MRD codes was given by Trombetti and Zhou in \cite{trombetti2018new}, and it was based on the same auxiliary result that we generalized in Proposition \ref{prop:full_kernel_norm_condition}. In this section we provide the same construction for the sum-rank metric, and show that they are as well MSRD. As a byproduct, this construction also produces a new family of additive MDS codes.\footnote{An additive code is a code which is linear over the prime field.}

\medskip

\subsection{Twisted Linearized Reed-Solomon Codes of TZ-Type}
Since we have that $\LL/\K$ is a cyclic Galois extension of degree $n$, then for every divisor $s$ of $n$ there exists an intermediate field   $\K\subseteq \E \subseteq \LL$ such that $[\LL:\E]=s$ and $\E/\K$ is Galois with $\Gal(\E/\K)=\langle \theta^s\rangle$.\footnote{Here, with a slight abuse of notation, we are writing $\theta^t$ using  $\theta$ also to denote the restriction of $\theta$ to $\E$.} Moreover, we will write $\K^{(2)}$ to denote the multiplicative subgroup of $\K^*$ consisting of all the squares, that is
$$\K^{(2)}:=\{a^2 \colon a \in \K^*\}.$$

\begin{definition}\label{def:TZ_codes}
 Let $n=2t$ and let $\E$ be an intermediate extension such that $[\LL:\E]=2$. Let $\gamma \in \LL^*$ be such that $\Norm(\gamma)\notin\K^{(2)}$. Moreover, assume that $\Lambda\subseteq\K^{(2)}$. The code
 $$\mD_k^\theta(\gamma):=\left\{ f_0+ \ldots+ f_{k-1}X^{k-1} +\gamma f_k X^{k} \colon f_1,\ldots,f_{k-1} \in \LL, f_0,f_k \in \E \right\}\subseteq \qspace $$
 is called \textbf{twisted linearized Reed-Solomon code of TZ-type}.
\end{definition}

Also in this case we can prove that these codes are MSRD using Proposition \ref{prop:full_kernel_norm_condition}.

\begin{theorem}\label{thm:TZ_are_MSR}
 The code $\mD_k^{\theta}(\gamma)$ is a maximum sum-rank distance code.
\end{theorem}

\begin{proof}
First, observe that the code $\mD_k^{\theta}(\gamma)$ is $\E$-linear, and $\dim_{\E}(\mD_k^{\theta}(\gamma))=2k$. Hence, by Theorem \ref{thm:singleton}, we need to show that $\dd_{\Lambda}(\mD_k^{\theta}(\gamma))=\ell n-k+1$.
Let $F(X)=f_0+ \ldots+ f_{k-1}X^{k-1} +\gamma f_k X^{k}\in \mD_k^{\theta}(\gamma)$ be a nonzero skew polynomial. If $f_k=0$, then $\deg F(X)\leq k-1$ and by Theorem \ref{thm:bound_degree_improved} we have $\wt_{\Lambda}(F)\geq \ell n-k+1$. Hence, assume that $f_k\neq 0$. If by contradiction we have $\wt_{\Lambda}(F)=\ell n-k$, this means that $\deg(F)=\dim_{\K}(\ker(F))=k$. By Proposition \ref{prop:full_kernel_norm_condition} we must have
$$ \Norm(f_0/\gamma f_k)=(-1)^{2tk}\prod_{i=1}^\ell \lambda_i^{d_{\lambda_i}(F)}=\prod_{i=1}^\ell \lambda_i^{d_{\lambda_i}(F)} \in \K^{(2)}. $$
On the other hand, we have
\begin{align*}\Norm(f_0/\gamma f_k)&=\Norm(\gamma)^{-1} \Norm(f_0/f_k)\\
&=\Norm(\gamma)^{-1}\N_{\E/\K}(\N_{\LL/\E}(f_0/f_k))\\
&=\Norm(\gamma)^{-1}\N_{\E/\K}(f_0/f_k)^2 \notin \K^{(2)},
\end{align*}
which yields a contradiction. 
\end{proof}

\begin{remark}
 Also in this case, if we fix $\ell=1$ and choose $\bs\alpha=\alpha_1=1$, the resulting code coincides with the code introduced by Trombetti and Zhou in \cite{trombetti2018new}, with respect to the rank metric. 
\end{remark}

We now state the results on the dual and the adjoint of a  twisted linearized Reed-Solomon code of TZ-type. The proof are omitted since they can be derived with straightforward computations as done in Theorems \ref{thm:dual_twisted} and \ref{thm:adjoint_twisted}.

\begin{proposition}
The dual code of a twisted linearized Reed-Solomon code of TZ-type is given by
 $$\mD_k^\theta(\gamma)^\perp = \mD_{\ell n-k}^\theta(-\gamma)\cdot X^k.$$
\end{proposition}

\begin{proposition}
Let $\Lambda$ be a cyclic group. Then, 
the adjoint code of a twisted linearized Reed-Solomon code of TZ-type is given by
$$\mD_k^\theta(\gamma)^\top=\mD_{k}^{\theta}(\theta^{n-k}(\gamma^{-1}))\cdot X^{\ell n-k}.$$
\end{proposition}

As we did for twisted linearized Reed-Solomon codes in Proposition \ref{prop:number_of_blocks_twisted}, here we determine the possible values $\ell$ and the existence of  twisted  linearized Reed-Solomon codes of TZ-type.

\begin{proposition}\label{prop:number_of_blocks_TZtwisted}\,
\begin{enumerate}
    \item If $\K=\Fq$, then we can construct a  twisted linearized Reed-Solomon code of TZ-type $\mD_k^\theta(\gamma)$ with $\ell$ blocks whenever $q$ is odd and  $\ell \leq \frac{q-1}{2}$.
    \item If $\K$ is an infinite field with $\mathrm{char}(\K)\neq 2$ admitting a cyclic Galois extension $\LL$ of degree $n$, then we can construct a  twisted linearized Reed-Solomon code of TZ-type $\mD_k^\theta(\gamma)$ with $\ell$ blocks, for every $\ell \in \mathbb N$.
    \end{enumerate}
\end{proposition}

\begin{proof}
 By definition, for the existence of an element $\gamma$ such that $\Norm(\gamma) \notin \K^{(2)}$, we need that $\K^{(2)} \subsetneq \K^*$, and this is possible if and only if $\mathrm{char}(\K) \neq 2$.  
 \begin{enumerate}
     \item  In the case of finite fields, the  subgroup $\Fq^{(2)}$ has cardinality $\frac{q-1}{2}$, and therefore we can choose as $\Lambda$ a set of $\ell$ elements for any  $\ell \leq \frac{q-1}{2}$.
     \item If $\K$ is infinite and $\mathrm{char}(\K)\neq 2$, then $|\K^*/\K^{(2)}|=2$, and thus we can choose as $\Lambda$ any set of $\ell$ elements in $\K^*$, for every $\ell \in \mathbb N$. 
 \end{enumerate}
\end{proof}

\medskip

\subsection{New Additive MDS Codes}
The importance of MDS codes in coding theory is undeniable. These are codes meeting the Singleton bound with equality, and they have the greatest error correction and detection capabilities in the Hamming metric, for the given size and length. The most prominent faily of MDS codes is certainly given by Reed-Solomon codes. These codes are linear over the underlying field, as many of the codes studied in the literature, due to the efficiency of encoding and decoding operations. However, also additive MDS codes also play an important role and they have been investigated also for their geometric equivalent interpretation as arcs of projective subspaces; see \cite{ball2020additive}.
In addition, it was shown in  \cite{ketkar2006nonbinary} that the existence of additive MDS codes is equivalent to the existence of quantum stabilizer MDS codes. 

We now focus on the special case of our code framework in which the block dimension is $n=1$. Although it can be derived from the general setting explained in Definition \ref{def:TZ_codes}, here we give a detailed description, since what we are going to obtain is a new family of additive MDS codes and we want to make our best to have a self-contained subsection. We will  focus only on the finite field case, since MDS codes over infinite fields are probably less interesting.

The special case of block dimension equal to $1$ coincides with the Hamming metric. The connection is easily explained by specializing the isomorphism given in Theorem \ref{thm:isomorphism_skew_sumrank} to the case of commutative polynomials with $\K=\LL$ and $\theta=\mathrm{id}$. Let $\Lambda=\{\lambda_1,\ldots,\lambda_\ell\}\subseteq \LL$ be a set of distinct evaluation points, and let us denote by $\dd_{\HH}$ the \emph{Hamming distance} on $\LL^\ell$, that is the distance map defined as
$ \dd_{\HH}(\mathbf{u},\mathbf{v})=\wt_{\HH}(\mathbf{u}-\mathbf{v})$, where $$\wt_{\HH}((v_1,\ldots,v_{\ell}))=|\{i \colon v_i \neq 0 \}|.$$
In this case, we have that 
$$H_{\Lambda}(X)=\prod_{i=1}^\ell(X-\lambda_i),$$
and the isomorphism of Theorem \ref{thm:isomorphism_skew_sumrank} can be rewritten as 
$$\begin{array}{rccl} \ev_{\bs\lambda}:&\LL[X]/(H_{\Lambda}(X)) &\longrightarrow &\LL^\ell, \\
& F(X) & \longmapsto & (F(\lambda_1),\ldots,F(\lambda_\ell)),\end{array} $$
where $\bs\lambda=(\lambda_1,\ldots,\lambda_\ell)$
This isomorphism is also an isometry, since one can verify, using the Chinese Remainder Theorem, that
$$\wt_{\Lambda}(F(X))=|\{i \colon F(\lambda_i)\neq 0\}|=\wt_{\HH}(\ev_{\bs\lambda}(F(X))).$$

We now fix $q$ to be an odd prime power and we assume now that $\LL=\K=\F_{q^2}$.
Moreover, we also consider the subfield $\E=\Fq$. The only difference with respect to Definition \ref{def:TZ_codes} is that in this case here the field $\E=\Fq$ is strictly contained in $\K=\F_{q^2}$ and not viceversa. Nevertheless, this is not going to affect the construction, which is given below.

\begin{definition}
 Let $1\leq k \leq \ell-1$ be positive integers. Let $\gamma \in \F_{q^2}\setminus(\F_{q^2})^{(2)}$, and let $\Lambda \subseteq (\F_{q^2})^{(2)}$ with $|\Lambda|=\ell$. The code 
 $$\mD_k(\gamma):=\left\{ f_0+ \ldots+ f_{k-1}X^{k-1} +\gamma f_k X^{k} \colon f_1,\ldots,f_{k-1} \in \F_{q^2}, f_0,f_k \in \Fq \right\}\subseteq \F_{q^2}[X]/(H_{\Lambda}(X)) $$
 is called a \textbf{Twisted Reed-Solomon code of TZ-type}.
\end{definition}

We now prove that these codes are maximum distance separable (MDS).

\begin{theorem}\label{thm:TZ_MDS} 
  The code $\ev_{\bs\lambda}(\mD_k(\gamma))$ is an $\Fq$-linear code of length $\ell$, size $q^{2k}$ and minimum Hamming distance $\ell-k+1$. In other words, $\ev_{\bs\lambda}(\mD_k(\gamma))$ is an $\Fq$-linear $(\ell,q^{2k},\ell-k+1)_{q^2}$ MDS code.
\end{theorem}

\begin{proof}
It is clear that the set $\mD_k(\gamma)$ is $\Fq$-linear and has size $q^{2k}$. Moreover, the map $\ev_{\bs\lambda}$ is clearly injective.  We only need to prove that, for each $F(X)=f_0+ \ldots+ f_{k-1}X^{k-1} +\gamma f_k X^{k}\in \mD_k(\gamma)$, we have $|\{i \colon F(\lambda_i)=0\}|\leq k+1$. If $f_k=0$, then $\deg F(X) \leq k-1$, and hence clearly  $|\{i \colon F(\lambda_i)=0\}|\leq k+1$. Now, assume that $f_k\neq 0$ and $|\{i \colon F(\lambda_i)=0\}|\leq k+1$. Therefore, there exists a subset $\Lambda'\subseteq \Lambda$ with $|\Lambda'|=k$ such that
$$ (\gamma f_k)^{-1}F(X)=\prod_{\lambda\in\Lambda'}(X-\lambda).$$
Thus, the degree-$0$ coefficient of the polynomial $(\gamma f_k)^{-1}F(X)$ is equal to 
$$(-1)^{k}\prod_{\lambda \in\Lambda'}\lambda \in (\F_{q^2})^{(2)}.$$
On the other hand, the degree-$0$ coefficient of the polynomial $(\gamma f_k)^{-1}F(X)$ is also equal to
$ f_0(f_k\gamma)^{-1}$
which does not belong to $(\F_{q^2})^{(2)}$, since $f_0/f_k\in \Fq^*$ and $\Fq^*\subseteq (\F_{q^2})^{(2)}$. This yields a contradiction and proves the claim. 
\end{proof}

\section{Conclusions and Future Work}\label{sec:conclusions}

In this paper we have proposed a new natural algebraic framework to study codes endowed with the sum-rank metric. The ambient space is a $\K$-algebra obtained as the quotient of a skew polynomial ring over a field extension $\LL$ of $\K$ by a suitable two-sided ideal. This space is then proved to be isometric to the classical vector framework endowed with the sum-metric, as well as the matrix framework. This result is based on recent works by McGuire and Sheekey \cite{mcguire2019characterization}; see also \cite{caruso2019residues}.  We have then studied sum-rank metric codes in this setting and their duality theory. Our approach appears very natural as a generalization of the polynomial approach for codes in the Hamming metric, and of the $q$-linearized polynomial approach for rank-metric codes over finite fields. We have then introduced \emph{twisted linearized Reed-Solomon codes}, which are the  counterpart of Sheekey's twisted Gabidulin codes. With an argument that resembles the proof that twisted Gabidulin codes are optimal in the rank metric (i.e. MRD), we were able to prove that twisted linearized Reed-Solomon codes are maximum sum-rank distance codes. The same strategy also allows to show that the same construction given by Trombetti and Zhou (see \cite{trombetti2018new}) in the rank metric setting produces here a second family of MSRD codes. As a byproduct, in this way we derive a new family of $\Fq$-linear  MDS codes in the Hamming metric  over $\F_{q^2}$  of length up to $\frac{q^2-1}2$.

In light of this new point of view, it is natural to ask which properties sum-rank metric codes share with rank-metric and Hamming-metric codes, and which constructions of codes can be borrowed. 
In \cite{longobardizanellascatt} and \cite{longobardi2021large}, two new families of MRD codes over finite fields were introduced, which were recently generalized to a new larger family in \cite{neri2021extending}. These codes, in the language of skew polynomials, can be described as follows. Let $n=2t$ with $t \geq 3$, $q$ be an odd prime power and fix any $h \in \F_{q^n}$ such that $\N_{q^n/q^t}(h)=-1$. Define the skew polynomial 
$$\psi_{h,t}(X)=X +X^{{t-1}}+h\theta(h)X^{t+1}+h\theta^{-1}(h^{-1})X^{{2t-1}}.$$
Then, the rank-metric code
$$  \langle 1, \psi_{h,t}(X) \rangle_{\F_{q^n}} \subseteq \Fn[X;\theta]/(X^n-1),$$
is an MRD code.  It would be interesting to determine whether one can adapt this construction to develop a more general family of MSRD codes in $\qspace$. 

\section*{Acknowledgements}

The author would like to thank Gianira N. Alfarano and {John Sheekey} for fruitful discussions and useful comments about the manuscript.

\bibliographystyle{abbrv}
\bibliography{biblio}

\end{document}